\documentclass[a4paper, reqno]{amsart}

\usepackage[utf8]{inputenc}
\usepackage[T1]{fontenc}
\usepackage{lmodern}
\usepackage[english]{babel}

\usepackage[usenames,dvipsnames,svgnames,table]{xcolor}
\usepackage{hyperref}

\usepackage{amssymb}

\usepackage{mathtools}
\usepackage{graphicx}
\usepackage{nicefrac}
\usepackage[symbol*]{footmisc}
\usepackage{tikz}
\usepackage{tikz-cd}

\usepackage{cleveref,nameref}
\usepackage{xspace}

\usepackage{enumitem}

\usepackage{xifthen}
\usepackage{verbatim}

\DeclareMathOperator{\psc}{psc}
\DeclareMathOperator{\torp}{torp}
\DeclareMathOperator{\rot}{rot}

\newenvironment{optional}%
        {%
            \ifthenelse{\isundefined{\alsoShowOptionalContent}}%
                    {\expandafter\comment}%
                    {\color{red}}%
                    }%
         {%
            \ifthenelse{\isundefined{\alsoShowOptionalContent}}%
                    {\expandafter\endcomment}%
                    {}%
          }

\crefname{section}{\S\!}{\S\S\!}
\crefname{satz}{Theorem}{Theorems}
\Crefname{satz}{Theorem}{Theorems}
\crefname{msatz}{Theorem}{Theorems}
\Crefname{msatz}{Theorem}{Theorems}
\crefname{lem}{Lemma}{Lemmata}
\Crefname{lem}{Lemma}{Lemmata}
\crefname{prop}{Proposition}{Propositions}
\Crefname{prop}{Proposition}{Propositions}
\crefname{kor}{Corollary}{Corollaries}
\Crefname{kor}{Corollary}{Corollaries}
\crefname{mkor}{Corollary}{Corollaries}
\Crefname{mkor}{Corollary}{Corollaries}
\crefname{bem}{Remark}{Remarks}
\Crefname{bem}{Remark}{Remarks}
\crefname{dfn}{Definition}{Definitions}
\Crefname{dfn}{Definition}{Definitions}
\crefname{bsp}{Example}{Examples}
\Crefname{bsp}{Example}{Examples}

\newcommand{\halfrule}{\par\parbox{\textwidth}{\noindent\hfil\rule{0.5\textwidth}{.4pt}\hfil}\par}

\newcommand{\blank}{{\,\cdot\,}} % blank argument

\iflanguage{ngerman}{
  \DeclareMathOperator{\tr}{spur}

}{

  \DeclareMathOperator{\tr}{tr}

}

\DeclareMathOperator{\End}{End} % endomorphisms
\DeclareMathOperator{\Emb}{Emb} % embeddings
\DeclareMathOperator{\Diff}{Diff}

\DeclareMathOperator{\Orthogonal}{O}

\DeclareMathOperator{\Spin}{Spin}

\DeclareMathOperator{\KOTheory}{KO}

\DeclareMathOperator{\Classifying}{B}

\DeclareMathOperator{\Ric}{Ric} % Ricci Tensor field
\DeclareMathOperator{\scal}{scal} % scalar curvature

\DeclareMathOperator{\id}{id}

\DeclareMathOperator{\pr}{pr} % canonical projection on a subscript component

\newcommand{\Ordnung}[1][]{ % Ordnung einer Gruppe
  \ifthenelse{\isempty{#1}}{
    \#
  }{
    \left|#1\right|
  }
}
\renewcommand{\Im}{\ensuremath{\operatorname{Im}}}
\DeclareMathOperator{\R}{\mathbb{R}}
\DeclareMathOperator{\E}{\mathbb{E}}
\DeclareMathOperator{\N}{\mathbb{N}}

\DeclareMathOperator{\Tangent}{T} % tangent bundle
\makeatletter
\providecommand*{\dop}%
{\@ifnextchar^{\@diffOperator}{\@diffOperator^{}}}
\def\@diffOperator^#1{%
  \mathop{\mathrm{\mathstrut d}}%
  \nolimits^{#1}\gobblespace}

\providecommand*{\Dop}%
{\@ifnextchar^{\@DiffOperator}{\@DiffOperator^{}}}
\def\@DiffOperator^#1{%
  \mathop{\mathrm{\mathstrut D}}%
  \nolimits^{#1}\gobblespace}

\def\gobblespace{%
  \futurelet\diffarg\opspace}
\def\opspace{%
  \let\DiffSpace\!%
  \ifx\diffarg(%
  \let\DiffSpace\relax
  \else
  \ifx\diffarg[%
  \let\DiffSpace\relax
  \else
  \ifx\diffarg\{%
  \let\DiffSpace\relax
  \fi\fi\fi\DiffSpace}
\makeatother

\let\emptyset\varnothing

\newcounter{thmmain}

\numberwithin{equation}{section}

\theoremstyle{plain}
\newtheorem{satz}[equation]{Theorem}
\newtheorem{msatz}[thmmain]{Theorem}
\newtheorem{lem}[equation]{Lemma}
\newtheorem{kor}[equation]{Corollary}
\newtheorem{mkor}[thmmain]{Corollary}
\newtheorem{prop}[equation]{Proposition}

\theoremstyle{definition}
\newtheorem{dfn}[equation]{Definition}
\newtheorem{bsp}[equation]{Example}

\theoremstyle{remark}
\newtheorem{bem}[equation]{Remark}

\theoremstyle{plain}
\newtheorem*{satz*}{Theorem}
\newtheorem*{lem*}{Lemma}
\newtheorem*{kor*}{Corollary}
\newtheorem*{prop*}{Proposition}
\newtheorem*{conv*}{Convention}

\theoremstyle{definition}
\newtheorem*{dfn*}{Definition}
\newtheorem*{bsp*}{Example}

\theoremstyle{remark}
\newtheorem*{bem*}{Remark}

\swapnumbers

\newtheorem*{acknowledgements}{Acknowledgements}

\setlist[enumerate]{label=(\arabic*), itemsep=0cm, leftmargin=1cm}
\setlist[itemize]{itemsep=0cm, leftmargin=1cm}

\hypersetup{
  pdfcreator = {},
  pdfproducer = {}
}
\usepackage[explicit]{titlesec}
\titleformat{\section}{\normalfont\large\filcenter}{\S \thesection.}{3.33pt}{#1}
\titleformat{name=\section,numberless}{\normalfont\large\filcenter}{}{0cm}{#1}
\titlespacing{\section}{0cm}{0.4cm}{0.4cm}

\titleformat{\subsection}[runin]{\normalfont}{\S \textbf{\thesubsection}}{3.33pt}{\bfseries #1.}
\titleformat{name=\subsection,numberless}[runin]{\normalfont\bfseries}{}{3.33pt}{#1}
\titlespacing{\subsection}{0cm}{\parskip}{3.33pt}

\usepackage{titletoc}
\usepackage[titles]{tocloft}
\usepackage{calc}

\cftsetindents{section}{0cm}{0.7cm}

\usepackage{todonotes}

\tikzstyle{notestyleraw}=[
draw=\@todonotes@currentbordercolor,
fill=\@todonotes@currentbackgroundcolor,
line width=0.6pt,
text width=\@todonotes@textwidth-1.6ex-1pt,
inner sep=0.8ex
]

\colorlet{draftnotesbgcolor}{orange}
\presetkeys{todonotes}{linecolor=draftnotesbgcolor, backgroundcolor=draftnotesbgcolor, bordercolor=draftnotesbgcolor, figcolor=white}{}

\colorlet{draftsimplenotesbgcolor}{lightgray}

\title[On the space of metrics with surgery stable curvature]{On the space of
  riemannian metrics satisfying surgery stable curvature conditions}

\author[J.-B.~Kordaß]{Jan-Bernhard~Kordaß$^\star$}
\email{\href{mailto:jb@kordass.eu}{jb@kordass.eu}}
\thanks{${}^{\star}$This work was supported in part by the SNSF-Project 200021E-172469 and the DFG-Priority programme \emph{Geometry at infinity} (SPP 2026).
  The author acknowledges funding by the Deutsche Forschungsgemeinschaft (DFG, German Research Foundation) – 281869850 (RTG 2229).}
\urladdr{\url{http://www.kordass.eu}}
\date{September 15, 2020}

\begin{document}

\begin{abstract}
  We utilize a condition for algebraic curvature operators called \textit{surgery stability} as suggested by the work of S.\ Hoelzel to investigate the space of riemannian metrics over closed manifolds satisfying these conditions.
  Our main result is a parametrized Gromov-Lawson construction with not necessarily trivial normal bundles and shows that the homotopy type of this space of metrics is invariant under surgeries of a suitable codimension.
  This is a generalization of a well-known theorem by V.\ Chernysh and M.\ Walsh for metrics of positive scalar curvature.

  As an application of our method, we show that the space of metrics of positive scalar curvature on quaternionic projective spaces are homotopy equivalent to that on spheres.
\end{abstract}

\maketitle

\vspace{-0.5cm}

\section{Introduction}

\begin{optional}
  \makeatletter
  \providecommand\@dotsep{5}
  \def\listtodoname{List of Todos}
  \def\listoftodos{\@starttoc{tdo}\listtodoname}
  \makeatother

  \listoftodos
\end{optional}

In this work we aim to give a sufficient criterion for a curvature condition such that the space of such metrics is weakly homotopy equivalent to a subspace of prescribed form around a fixed submanifold and derive various conclusions.

A curvature condition is an open subset of the space of $\Orthogonal(n)$-invariant operators $\bigwedge^2 \mathbb E^n \to \bigwedge^2 \mathbb E^n$ satisfying the Bianchi identity, where $\mathbb E^n$ denotes the euclidean space.
% Based on the work of S.\ Hoelzel in \cite{hoelzel-2016}, we define conditions for a riemannian metric, which can be preserved under surgeries of codimension at least $c$, where $c \geq 3$ depends on the condition.
We will define two notions that we call \emph{deformable} and \emph{codimension $c$ surgery stable}, which allow to extend the Gromov-Lawson construction for curvature conditions developed by S.\ Hoelzel in \cite{hoelzel-2016} and will be explained in greater detail in \cref{subsec:surgery-stable,subsec:rotationally-symmetric-metrics}.

Let $M^n$ be a closed, smooth manifold of dimension $n$ and let $N^k$ be a compact, $k$-dimensional submanifold in $M$.
We denote by $\mathcal R_C(M)$ the space of riemannian metrics that satisfy the condition $C$, i.e.\ for any $p \in M$ and any linear isometry $i \colon \mathbb E^n \to (\Tangent_pM, g_p)$ we have $i^{*} R_p \in C$, where $R_p$ is the riemannian curvature operator of the riemannian metric $g$ at $p$.
Given a tubular map $\phi \colon \nu N \to M$ and a distinct metric on $\nu N$, we will define a space of metrics $R^{\torp}_C(M)$, which are \emph{standard near $N$}.

\begin{msatz}[parametrized Gromov-Lawson construction]\label{msatz:main-technical-introduction}
  Let $C \subset \mathcal C_{\mathrm B}(\E^n)$ be a deformable, codimension $c$ surgery stable curvature condition with $n - k \geq c$.
  Then the inclusion of metrics, which are standard near $N$
  \begin{align*}
    \mathcal R^{\torp}_C(M) \hookrightarrow \mathcal R_C(M)
  \end{align*}
  is a weak homotopy equivalence.
\end{msatz}

The condition positive scalar curvature, i.e.\ $C = (\psc)$, is deformable, codimension $3$ surgery stable and thus \Cref{msatz:main-technical-introduction} generalizes a result of V.\ Chernysh \cite{chernysh-2004} and M.\ Walsh \cite{walsh-2013} for submanifolds with trivial normal bundle, which has recently been revisited by J.\ Ebert and G.\ Frenck \cite{ebert-frenck-2018}.
Note that the normal bundle of $N$ does not necessarily have to be trivial in \Cref{msatz:main-technical-introduction}.
We remark that the Gromov-Lawson constructions given in \cite{schoen-yau-1979} and \cite{hoelzel-2016} also do not require a trivial normal bundle.

In the case of surgeries along embedded spheres with trivial normal bundles we immediately obtain from \cref{msatz:main-technical-introduction}:

\begin{msatz}\label{msatz:homotopy-invariance}
  Let $M_0$ be a closed $n$-manifold and let $M_1$ be obtained from $M_0$ by surgery of codimension $n-k \geq c$ with $k \geq c-1$.
  Then $\mathcal R_C(M_0)$ is homotopy equivalent to $\mathcal R_C(M_1)$.
\end{msatz}

\cref{msatz:main-technical-introduction} can be used to consider more general gluing constructions, namely surgeries along embedded spheres with nontrivial normal bundles.
In particular, we can cut along the sphere bundle of the Hopf fibration $S^{4k+3} \to \mathbb HP^k$ within $\mathbb HP^{k+1}$ to obtain $S^{4(k+1)}$ from a generalized surgery of codimension $4$.
\begin{msatz}\label{msatz:hpn}
  Let $k \geq 1$.
  Then
  \begin{align*}
    \mathcal R_{\psc}(\mathbb HP^k) \simeq \mathcal R_{\psc}(S^{4k})
    \quad\text{and}\quad
    \mathcal R_{1\text{-curv}>0}(\mathbb HP^k) \simeq \mathcal R_{1\text{-curv}>0}(S^{4k}).
  \end{align*}
\end{msatz}

By an application of the original Gromov-Lawson construction (\cite{gromov-lawson-1980} and independently \cite{schoen-yau-1979}) we know that a closed, simply-connected manifold of dimension $n \geq 5$, which is oriented bordant (resp. spin bordant, if $M$ is spin) to a manifold with positive scalar curvature admits a metric of positive scalar curvature itself.
The proof of this result relies on the fact that one can arrange the bordism to be composed from traces of codimension $\geq 3$ surgeries, which allow to push through positive scalar curvature.
It is well-known that this procedure can be generalized to other tangential structures.
Denoting by $\Classifying \Orthogonal\left<l\right>$ the $l\,$th stage in the Whitehead tower of $\Classifying \Orthogonal$, a closed manifold $M^n$ is called \emph{$\Classifying \Orthogonal\left<l\right>$-manifold}, if the stable normal bundle $M \to \Classifying \Orthogonal$ lifts to $M \to \Classifying \Orthogonal\left<l\right>$ through the tower.
The corresponding bordism group will be denoted by $\Omega^{\left<l\right>}_n$.

\begin{prop*}\footnote{e.g.\ \cite[Proposition 3.4]{botvinnik-labbi-2014}}
  Let $M_0, M_1$ be closed $\Classifying \Orthogonal\left<l\right>$-bordant $n$-manifolds and let $r \geq 1$ be such that $n \geq 2r + 3$ and $l \geq r + 2$.
  If $M_1$ is $r$-connected, then $M_1$ can be obtained from $M_0$ by surgeries of codimension at least $r + 2$.
\end{prop*}

Combined with the Gromov-Lawson construction for curvature conditions \cite[Theorem A]{hoelzel-2016} we obtain:

\begin{msatz}\label{msatz:main}
  Let $(M_0,g)$ be a closed riemannian manifold satisfying a codimension $c$ surgery stable curvature condition, which is $\Classifying \Orthogonal\left<l\right>$-bordant to an $r$-connected, closed manifold $M_1$.
  If $n \geq 2r + 3$, $l \geq r + 2$ and $r + 2 \geq c$, then $M_1$ admits a riemannian metric satisfying the same curvature condition as $M_0$.
\end{msatz}

Moreover, using \cref{msatz:main-technical-introduction}, we can strengthen this statement for highly-connected manifolds.

\begin{msatz}
  Let $C$ be a deformable,  codimension $c$ surgery stable curvature condition.
  Let $M_0, M_1$ be $\Classifying \Orthogonal\left<l\right>$-bordant, $r$-connected, closed manifolds with $n \geq 2r + 3$, $l \geq r + 2$ and $r + 2 \geq c$.
  Then $\mathcal R_C(M_0)$ is homotopy equivalent to $\mathcal R_C(M_1)$.
\end{msatz}

Let us state the result for the following examplary case.
A riemannian manifold $(M, g)$ is is said to have \emph{positive $1$-curvature}%
\footnote{{%
  This condition is sometimes also called \emph{positive Einstein} or the metric is said to have \emph{positive Einstein tensor}.%
  We will refrain from these terms to avoid confusion with \emph{Einstein} metrics with positive Einstein constant.}}, if $\scal_g - 2\Ric_g(v) > 0$ for all points $p \in M$ and unit vectors $v \in \Tangent_pM$.
It is not hard to see that every metric of positive $1$-curvature has positive scalar curvature.
We will see later that positive $1$-curvature defines a codimension $4$ surgery stable curvature condition $C \coloneqq (1\text{-curv}>0)$.
In particular, we can consider the space of metrics with positive $1$-curvature $\mathcal R_{1\text{-curv}>0}(M)$ over a fixed manifold $M$.
Noting that $\Classifying \Orthogonal\left<4\right> = \Classifying \Spin$ and the spin bordism group in dimension $7$ is trivial, we conclude the following.

\begin{mkor}
  For every $2$-connected, $7$-dimensional, closed spin manifold $M$ the space of riemannian metrics with positive $1$-curvature $\mathcal R_{1\text{-curv}>0}(M)$ has the homotopy type of $\mathcal R_{1\text{-curv}>0}(S^7)$.
\end{mkor}

This is quite a large class of manifolds (cf.\ the classification by topological invariants \cite{crowley-nordstroem-2018}), but in particular it holds for all exotic seven spheres.
Typical interesting examples also include $S^3$-bundles over $S^4$, which contain a family of infinitely many mutually distinct homotopy types admitting metrics of non-negative sectional and positive Ricci curvature (cf.\ \cite[Proposition 3.3]{grove-ziller-2000}).

It is a well-known consequence of the Atiyah-Singer index theorem that a homotopy sphere is spin-nullbordant if its $\alpha$-invariant vanishes (cf.\ \cite[p.144]{lawson-michelsohn-1989}).

\begin{mkor}
  Let $\Sigma^n$ be a homotopy sphere with $\alpha(\Sigma^n) = 0$ (e.g.\ in the case $n \not\equiv 1,2$ mod $8$) and with $n \geq 7$.
  Then the space of riemannian metrics with positive $1$-curvature $\mathcal R_{1\text{-curv}>0}(\Sigma)$ has the homotopy type of $\mathcal R_{1\text{-curv}>0}(S^n)$.
\end{mkor}

The existence question for metrics of positive $1$-curvature on $2$-connected manifolds was already addressed in \cite[Theorem I]{labbi-1997} and the homotopy equivalence of $\mathcal R_{\psc}(\Sigma^n) \simeq \mathcal R_{\psc}(S^n)$ is a corollary to the version of \cref{msatz:main} for positive scalar curvature ($\psc$) by Chernysh and Walsh (cf. \cite[Corollary 4.2]{walsh-2013}).
Note however that the inclusion $\mathcal R_{1\text{-curv}>0}(M) \hookrightarrow \mathcal R_{\psc}(M)$ is not understood.

The above corollaries hold in analogy for arbitrary curvature conditions, which are stable under surgery of codimension at least $4$.

\makeatletter
\providecommand\@dotsep{5}
\makeatother

\tableofcontents

\section{Preliminaries}

\subsection{Curvature conditions and spaces of riemannian metrics}

We will briefly recall the definition of algebraic curvature operators and notions of curvature.
These let us define what we understand under a curvature condition satisfied by a riemannian metric.

\begin{dfn}
  Let $(M,g)$ be a riemannian manifold with riemann curvature tensor $R \colon \mathcal V(M) \times \mathcal V(M) \to \End(\mathcal V(M))$, where $\mathcal V(M)$ denotes the space of smooth vector fields on $M$.
  One obtains what is called \emph{curvature operator} $\overline R \colon \bigwedge^2 \mathcal V(M) \to \bigwedge^2 \mathcal V(M)$ of $M$ via the relation
  \begin{align*}
    g(\overline R(X \wedge Y), Z \wedge W) = g(R(X,Y)W, Z),
  \end{align*}
  where $g$ on the left hand side is the extension of $g$ to $\Gamma(M, \bigwedge^2 \Tangent M)$.

  At every point $p \in M$ the curvature operator $\overline R$ defines a self-adjoint endomorphism
  \begin{align*}
    \textstyle \overline R_p \colon \bigwedge^2 \Tangent_pM \to \bigwedge^2 \Tangent_pM.
  \end{align*}
\end{dfn}

We denote by $\E^n$ the euclidean space $\R^n$ endowed with the standard inner product $\langle\blank,\blank\rangle$.
A self-adjoint endomorphism $R \colon \bigwedge^2 \E^n \to \bigwedge^2 \E^n$ is called \emph{algebraic curvature operator}.
From the inner product, we obtain the isomorphism
\begin{align*}
  \textstyle \eta \colon \bigwedge^2 \E^n \to \mathfrak{so}(n) \subset \End(\R^n),
  \quad
  x \wedge y \mapsto - \left<x,\blank\right>y + \left<y, \blank\right>x
\end{align*}
and thus for every $x,y \in \R^n$, an algebraic curvature operator gives rise to a skew-symmetric endomorphism
\begin{align*}
  R(x,y) \colon \R^n \to \R^n,
  \quad
  z \mapsto (\eta \circ R (x \wedge y))(z).
\end{align*}
An algebraic curvature operator $R$ is said to satisfy the \emph{Bianchi indentity}, if
\begin{align*}
  R(x, y)z + R(y, z)x + R(z, x)y = 0
\end{align*}
for all $x,y,z \in \R^n$.

Exactly in the same manner as we define sectional, Ricci and scalar curvature in a tangent space, we let
\begin{enumerate}
\item $\sec(R, E) \coloneqq \left< R(x,y)y, x \right>$ for an orthonormal basis $\{x,y\}$ of a 2-plane $E \leq \mathbb E^n$,
\item $\Ric(R, z) \coloneqq \sum_{i = 1}^n \left< R(e_i,z)z, e_i \right>$ for $z \in \mathbb E^n$ with $\|z\| = 1$,
\item $\scal(R) \coloneqq \sum_{i,j = 1}^n \left< R(e_i, e_j)e_j, e_i \right> = 2 \tr R$,
\end{enumerate}
where $\{e_i\}$ is the standard basis of $\mathbb E^n$.\footnote{%
  As common notation suggests, we will write relations such as $\sec(R) < \alpha$ to mean ``$\sec(R, E) < \alpha$ for every plane $E < \mathbb E^n$''.}

Let $(M,g)$ be a riemannian manifold and let $i \colon \E^n \to (\Tangent_pM, g_p)$ be a linear isometry into the tangent space at some point $p$ in $M$.
This defines an algebraic curvature operator
\begin{align*}\textstyle
  i^{*}\overline R_p \coloneqq (i^{\wedge 2})^{-1} \circ R_p \circ i^{\wedge 2} \colon \bigwedge^2 \E^n \to \bigwedge^2 \E^n
\end{align*}
for $i^{\wedge 2} \colon \bigwedge^2 \E^n \to \bigwedge^2 \Tangent_pM$ and $(i^{\wedge 2})^{-1} \colon \bigwedge^2 \Tangent_pM \to \bigwedge^2 \E^n$ induced by $i$ and its inverse.
Clearly, $i^{*}\overline R_p$ satisfies the Bianchi identity and we have $\sec(i^{*}\overline R_p, E) = \sec^g_p(i(E))$, $\Ric(i^{*}R_p,z) = \Ric^g_p(i(z))$ and $\scal(i^{*}R_p) = \scal^g(p)$.
Hence, we can use an algebraic curvature operator, which satisfies the Bianchi identity, to describe curvature properties of $g$ at $p$ up to the choice of an orthonormal basis in $\Tangent_pM$.

Let $\mathcal C_{\mathrm B}(\E^n)$ denote the vector space of algebraic curvature operators satisfying the Bianchi identity.
Then $\Orthogonal(n)$ acts on $\mathcal C_{\mathrm B}(\E^n)$ by change of the orthonormal basis in $\mathbb E^n$, that is via
\begin{align*}
  (A,  R) \mapsto (A \wedge A)^{-1} \circ R \circ (A \wedge A).
\end{align*}

\begin{dfn}
  A \emph{curvature condition} is an open subset $C$ of $\mathcal C_{\mathrm B}(\E^n)$, which is invariant under the action of $\Orthogonal(n)$ on $\mathcal C_{\mathrm B}(\E^n)$.
  We say that a riemannian metric $g$ on $M$ \emph{satisfies} $C$, if $i^{*}\overline R_p \in C$ for all linear isometries $i \colon \E^n \to (\Tangent_pM,g_p)$ and $p \in M$, where $\overline R_p$ denotes the curvature operator with respect to $g$ in $p$.
\end{dfn}

Thus, any curvature condition is a global statement on the structure of a riemannian manifold and obviously describes an isometry invariant property of a riemannian metric satisfying it.

\begin{bsp}
  \begin{enumerate}[label=(\roman*),leftmargin=1cm]
  \item Lower curvature bounds can be expressed as a curvature condition.
    For example, we can express (globally pointwise) positive sectional curvature as condition
    \begin{align*}
      (\sec > 0) \coloneqq \{ R \in \mathcal C_{\mathrm B}(\E^n) \mid \sec(R) > 0\},
    \end{align*}
    while positive scalar curvature is described by
    \begin{align*}
      \psc \coloneqq (\scal > 0) \coloneqq \{ R \in \mathcal C_{\mathrm B}(\E^n) \mid \scal(R) > 0 \}.
    \end{align*}
    We note that both sets are open cones in $\mathcal C_{\mathrm B}(\E^n)$ with $(\sec > 0) \subsetneq \psc$ for $n \geq 3$, while $(\sec > 0) = \psc$ for $n = 2$.

    \noindent
    We denote by $\mathbb S^n$ the $n$-sphere endowed with the round metric of radius $1$.
    For any isometry $i \colon \mathbb E^n \to \Tangent_p\mathbb S^n$, we have $i^{*}\overline R_p = \id_{\bigwedge^2\mathbb E^n}$.
    In particular $\mathbb S^n$ satisfies $(\sec > 0)$.
  \item Upper curvature bounds, such as
    \begin{align*}
      (\scal < \alpha) \coloneqq \{ R \in \mathcal C_{\mathrm B}(\E^n) \mid \scal(R) < \alpha\}
    \end{align*}
    for $\alpha \in \R$ are of course curvature conditions, although they will not play an important role in the following.
  \end{enumerate}
\end{bsp}

\begin{dfn}
  Let $M^n$ be a smooth manifold and denote by $\mathcal R(M)$ the space of complete riemannian metrics on $M$ equipped with the compact-open $C^{\infty}$-topology.
  Define, the \emph{space of riemannian metrics satisfying $C$} as
  \begin{align*}
    \mathcal R_C(M) \coloneqq \{ g \in \mathcal R(M) \mid g \text{ satisfies } C \} \subset \mathcal R(M).
  \end{align*}
  As curvature conditions are preserved under isometries, $\Diff(M)$ acts on $\mathcal R_C(M)$ by pullback of riemannian metrics, i.e.\ via
  \begin{align*}
    \Diff(M) \times \mathcal R_C(M) \to \mathcal R_C(M),
    \quad
    (\psi, g) \mapsto \psi^{*} g.
  \end{align*}
  The quotient $\mathcal M_C(M) \coloneqq \mathcal R_C(M)/\Diff(M)$ under this action is called \emph{moduli space of riemannian metrics satisfying $C$}.
\end{dfn}

\begin{bem*}
  These spaces are actually quite classical objects of interest in global riemannian geometry (cf. \cite[$\S\ 7\nicefrac{3}{4}$]{gromov-1991}).
  We refer to the survey \cite{tuschmann-wraith-2015} for further details and references regarding spaces of riemannian metrics.
\end{bem*}

Recently, there has been much progress showing that $\mathcal R_C(M)$ on certain closed manifolds is topologically non-trivial for $C = \psc$ and $C = (\Ric > 0)$ (cf. \cite{crowley-schick-2013}, \cite{hanke-schick-steimle-2014}, \cite{botvinnik-ebert-randal-williams-2017} and \cite{wraith-2011}).
For positive scalar curvature non-trivial homotopy elements can be exhibited factoring the Atiyah-Bott-Shapiro map of spectra $\mathcal A \colon \operatorname{M}\Spin \to \KOTheory$ through the space of metrics with positive scalar curvature carried out in \cite{botvinnik-ebert-randal-williams-2017}, which is a vast generalization of an idea by N.\ Hitchin.
Incidentally, already Hitchin's original construction (cf.  \cite[Section 4.4]{hitchin-1974}) identifies non-trivial homotopy elements in $\mathcal R_C(M)$ on any spin manifold $M$ of a certain dimension for every $C \subset \psc$ with $\mathcal R_C(M) \neq \emptyset$ that are produced using the orbit map of the above described $\Diff(M)$ action on $\mathcal R_C(M)$.

\begin{conv*}
  By a \emph{(compact) family of metrics} on a manifold $M$, we always understand the image of a continuous map $g \colon S \to \mathcal R(M)$ from a compact topological space $S$ and denote it by $\{g_{\xi}\}_{\xi \in S}$.
\end{conv*}

\subsection{Surgery stable curvature conditions}\label{subsec:surgery-stable}

\begin{dfn}
  Let $M$ be a smooth manifold of dimension $n$ and further let $\phi \in \Emb(S^k, M)$ be an embedding with trivial normal bundle.
  Thus $\phi$ extends to a tubular embedding $\phi \in \Emb(S^k \times D^{n-k}, M)$.
  Recall that a \emph{surgery} on $M$ along $\phi$ is a procedure, which endows the push-out
  \begin{equation*}
    \begin{tikzcd}
      S^k \times D^{n-k}(\frac{1}{2}, 1) \ar{r}{i} \ar[hook]{d}{\phi(\frac{1}{2}, 1)} & D^{k+1}(\frac{1}{2}, 1) \times S^{n-k-1} \ar{d}\\
      M \setminus \Im(S^k \times \mathring D^{n-k}(\frac{1}{2})) \ar{r} & \chi(M, \phi)
    \end{tikzcd}
  \end{equation*}
  where $\phi(\frac{1}{2}, 1) \coloneqq \phi|_{S^k \times D^{n-k}(\frac{1}{2}, 1)}$, $i \colon (\theta_1, \lambda \theta_2) \mapsto ((1-\lambda) \theta_1, \theta_2)$ and where we denote $D^k(a,b) \coloneqq \{ x \in D^k \mid a \leq |x| \leq b \}$,
  \begin{optional} i.e.
    \begin{align*}
      \chi(M, \phi) \coloneqq M \setminus \Im \mathring \phi \cup_{\partial \phi} D^{k+1} \times S^{n-k-1},
    \end{align*}
  \end{optional}
  with a differentiable structure compatible with the given differentiable structures on the individual parts.
  The number $n-k$ is called \emph{codimension of the surgery}.
\end{dfn}

Clearly, one can also equip $M$ with a riemannian metric and it is not hard to see that there exists a metric on $\chi(M,\phi)$, which coincides with the original metric away from the embedding $\phi$.
But since cutting and gluing of riemannian metrics along arbitrary submanifolds requires to smoothen the metric, it is not clear if this process keeps the curvature controlled.
It was only realized by Gromov and Lawson \cite{gromov-lawson-1980} and independently by R.\ Schoen and S.\ T.\ Yau \cite{schoen-yau-1979} that positive scalar curvature can be preserved under surgeries of codimension greater or equal three.

This theorem has been generalized in several directions and in particular to curvature conditions as defined above by Hoelzel.

\begin{dfn}\label{dfn:inner-cone-condition}
  A curvature condition $C \subset \mathcal C_{\mathrm B}(\E^n)$ is said to \emph{satisfy an inner cone condition with respect to} $S \in \mathcal C_{\mathrm B}(\E^n)\setminus\{0\}$, if there exists a continuous function $\rho \colon C \to (0, \infty)$ and for every $\rho = \rho(R)$ an open, convex $\Orthogonal(n)$-invariant cone $C_{\rho}$ containing $B_{\rho}(S)$ such that
  \begin{align*}
    R + C_{\rho} = \{ R + E \mid E \in C_{\rho} \} \subset C.
  \end{align*}
\end{dfn}

\begin{bem}\label{bem:curvature-condition-open-convex}
  \begin{enumerate}[label=(\roman*), leftmargin=1cm]
  \item \label{item:bem:curvature-condition-open-convex-i}
    If $C \subset \mathcal C_{\mathrm B}(\E^n)$ is a curvature condition given by an open, convex cone and $S \in C$, then $C$ satisfies an inner cone condition with respect to $S$.
  \item \label{item:bem:curvature-condition-open-convex-ii}
    If $C \neq \emptyset$ satisfies an inner cone condition with respect to $S \neq 0$, then there exists a $\lambda_0 > 0$ such that $\lambda S \in C$ for all $\lambda > \lambda_0$.

    \noindent The argument goes as follows.
    Let $R \in C$ be arbitrary and conclude from the inner cone condition that $R + C_{\rho} \in C$ for some open, convex $\Orthogonal(n)$-invariant cone with $B_{\rho}(S) \subset C_{\rho}$.
    Clearly, this implies $B_{\mu\rho}(R + \mu S) \subset C$ and the cone $\bigcup_{\mu > 0} B_{\mu\rho}(R + \mu S)$ intersects the line $\lambda S$ for some $\mu$ and $\lambda$ large enough.
  \end{enumerate}
\end{bem}

\begin{bem}
  Since the curvature operator of the standard sphere $S^q(1)$ of radius $1$, which we also denote by $\mathbb S^{q}$, is given by the identity, the curvature operator of $\E^{n-q} \times \mathbb S^q$ with the canonical product metric is precisely the projection map $R_{\E^{n-q} \times \mathbb S^q} \coloneqq \pi_{\bigwedge^2 \E^q} \colon \bigwedge^2\E^n \to \bigwedge^2\E^n$ induced by the projection $\E^n = \E^{n-q} \times \E^q \to \E^q$ on the last $q$ coordinates.
\end{bem}

\begin{dfn}
  Let $C \subset \mathcal C_{\mathrm B}(\E^n)$ be a curvature condition satisfying an inner cone condition with respect to $R_{\mathbb E^{n-c+1} \times \mathbb S^{c-1}}$ for some $c \in \{3, \ldots, n\}$.
  Then $C$ is said to \emph{admit codimension $c$ surgeries}.
\end{dfn}

\begin{prop}[{\cite[Proposition 2.2]{hoelzel-2016}}]
  If $C \subset \mathcal C_{\mathrm B}(\E^n)$ is a curvature condition satisfying an inner cone condition with respect to $R_{\mathbb E^{n - c + 1} \times \mathbb S^{c - 1}}$ for $3 \leq c \leq n$, then $C$ satisfies an inner cone condition with respect to $R_{\mathbb E^{n-c} \times \mathbb S^c}$.
\end{prop}

\begin{kor}
  Let $C \subset \mathcal C_{\mathrm B}(\E^n)$ be a curvature condition admitting codimension $c$ surgeries for some $c \in \{3, \ldots, n-1\}$.
  Then $C$ admits codimension $c+1$ surgeries.
\end{kor}

\begin{kor}\label{kor:round-sphere-satisfies-C}
  If $C$ admits codimension $c$ surgeries, then $C$ satisfies an inner cone condition with respect to $R_{\mathbb E^{n-q} \times \mathbb S^q}$ for all $c-1 \leq q \leq n$ and by \cref{bem:curvature-condition-open-convex} \ref{item:bem:curvature-condition-open-convex-ii} there exists a $\lambda_0 > 0$ such that $\lambda R_{\mathbb E^{n-q} \times \mathbb S^q} \in C$ for all $\lambda > \lambda_0$.
  In particular, $S^n(\frac{1}{\sqrt{\lambda}})$ with the round metric satisfies $C$.
\end{kor}

As the name suggests, admittance of surgery is precisely the assumption needed to obtain a new metric again satisfying this curvature condition on a surgery result.

\begin{satz}[{\cite[Theorem A]{hoelzel-2016}}]\label{satz:hoelzels-GL-construction}
  Let $C$ be a curvature condition admitting codimension $c$ surgeries and let $(M,g)$ be a riemannian manifold with $g$ satisfying $C$.
  If $\chi(M,\phi)$ is obtained from $M$ by surgery of codimension $\geq c$, then $\chi(M,\phi)$ admits a metric satisfying $C$.
\end{satz}

Thus we are only interested in the smallest value for which $C$ admits surgeries of such codimension leading us to the following definition.

\begin{dfn}
  A curvature condition $C$ admitting codimension $\tilde c$ surgeries is said to be \emph{codimension $c$ surgery stable}, if $c$ is minimal among all $\tilde c$.
\end{dfn}

\begin{optional}
  It might be conceptually better to say that $C$ is ``stable under codimension $c$ surgeries using the Gromov-Lawson construction'' or ``GL-surgery stable'', as there are curvature conditions that can be preserved by different means.
\end{optional}

\begin{bsp}\label{bsp:ss-cc-psc}
  \begin{enumerate}[label=(\roman*),leftmargin=1cm]
  \item By \cref{bem:curvature-condition-open-convex} \ref{item:bem:curvature-condition-open-convex-i}, a curvature condition $C$ given by an open, convex cone containing $R_{\mathbb E^{n-q+1} \times \mathbb S^{q-1}}$ for all $q \geq c$ (i.e.\ $\mathbb E^{n-q+1} \times \mathbb S^{q-1}$ with the standard product metric satisfies $C$) admits codimension $c$ surgeries.
  \item Positive scalar curvature $C = \psc$ is codimension $3$ surgery stable in this sense.
    Clearly, $c$ is minimal in the allowed range $\{3, \ldots n\}$ for $c$, but this conceptually makes sense as the standard metric on $\mathbb S^{2-1} \times \mathbb E^{n-2+1} = \mathbb S^1 \times \mathbb E^{n-1}$ is flat, i.e.\ does not have positive scalar curvature.
  \end{enumerate}
\end{bsp}

\begin{bsp}\label{bsp:ss-cc-p-curv}
  Another interesting curvature condition is \emph{positive $p$-curvature}, which interpolates between positive scalar (for $p = 0$) and positive sectional curvature (for $p = n-2$).
  It has been proposed by Gromov and was studied extensively by M.-L.\ Labbi (cf. \cite{labbi-1995,labbi-1997,labbi-1997-actions,labbi-2006}).

  It can be defined for $0 \leq p \leq n-2$ as an open convex cone
  \begin{align*}
    (p\text{-curv} > 0) \coloneqq \{ R \in \mathcal C_{\mathrm B}(\E^n) & \mid s_p(R)(P) > 0 \\
                                                                        & \quad \forall P \leq \R^n \text{ with } \dim P = p\},
  \end{align*}
  where $s_p(R) \colon G_p(\R^n) \to \R$ is the map $P \mapsto \sum_{i,j=1}^{n-p}\sec(R)(E_i,E_j)$ for an orthonormal basis $\{E_i\}$ of $P^{\perp}$ and $G_p(\R^n)$ is the real $p$-Grassmannian.

  By definition, positive $p$-curvature implies positive $(p-1)$-curvature and for every fixed dimension $n$ there is a sequence of cones
  \begin{align*}
    (\psc) = (0\text{-curv}) \supset (1\text{-curv}) \supset \cdots \supset ((n-3)\text{-curv}) \supset ((n-2)\text{-curv}) = (\sec > 0).
  \end{align*}
  Labbi showed in \cite{labbi-1997} that positive $p$-curvature is preserved under codimension $p+3$ surgeries.
  This is recovered by Hoelzel's theorem, as the condition is codimension $p+3$ surgery stable.

  Riemannian metrics satisfying this condition for values of $p$, which are small relative to $n$, exist in abundance.
  The product of every compact manifold with a $(p+3)$-sphere admits a metric of positive $p$-curvature.
  More specifically, every compact, connected Lie group with a bi-invariant metric (which is not the torus), every Einstein manifold with positive Einstein constant, as well as every Kähler manifold with positive Ricci curvature has positive $1$-curvature.
  On the other hand, there exist examples, which have positive $1$-curvature, while they do not admit any metric of positive Ricci curvature.
  Botvinnik and Labbi in \cite{botvinnik-labbi-2014} investigate obstructions to positive $p$-curvature for $p = 2,3$.
  For example, they show that a 3-connected non-string manifold $M$ of dimension at least 9 admit a metric with positive 2-curvature if and only if Hitchin's $\operatorname{KO}$-theoretic $\alpha$-invariant vanishes for $M$.
  Moreover, they find dimensions in which every 3-connected string manifold admits a metric with positive 2-curvature.
  Note however, that the situation is less clear in lower connectivity.
\end{bsp}

\begin{bsp}
  Similarly, one can consider a curvature condition, which interpolates between positive scalar (for $k = n$) and positive Ricci curvature (for $k = 1$) called \emph{$k$-positive Ricci curvature}.
  It was introduced by J.\ Wolfson in \cite{wolfson-2009} and can be defined for $1 \leq k \leq n$ as an open convex cone
  \begin{align*}
    (k\text{-pos Ric}) \coloneqq \{ R \in \mathcal C_{\mathrm B}(\E^n) & \mid \ \sum_{i=1}^k \Ric(R; e_i) > 0 \\
                                                                       & \quad \forall \{e_1, \cdots, e_k\} \text{ orthonormal} \}.
  \end{align*}
  As in the previous examples, there are successive inclusions $(k\text{-pos Ric}) \subset ((k+1)\text{-pos Ric})$ and in particular every metric with $k$-positive Ricci curvature for some $k$ has positive scalar curvature.
  Moreover, Wolfson showed that $k$-positive Ricci curvature for $2 \leq k \leq n-1$ is preserved under codimension $n-k+2$ surgeries.
  Surprisingly, $n$-positive Ricci curvature (positive scalar curvature), as well as $(n-1)$-positive Ricci curvature are both stable under surgeries of codimension $3$.
  It is an open question of Wolfson, if there exists a riemannian manifold of positive scalar curvature, which does not admit a metric of $(n-1)$-positive Ricci curvature.
\end{bsp}

Both of the curvature conditions decribed in the above examples can be regarded as \emph{intermediate curvature notions}, by which one might hope to understand the differences between both extremes in greater detail.

\subsection{Torpedo metrics}\label{subsec:rotationally-symmetric-metrics}

It is well-known that one can use warped products to describe rotationally symmetric metrics (cf. \cite[p.18ff]{petersen-2016}).
To do so consider a smooth function $\beta \colon [0,\infty) \to [0,\infty)$ and endow $(0,\infty) \times S^{q-1}$ with the metric $g^{\beta} \coloneqq \dop r^2 + \beta^2(r) g_{\mathbb S^{q-1}}$, where $g_{\mathbb S^{q-1}}$ is the round metric on the $q-1$ sphere.
If we assume that
\begin{enumerate}[label=(\roman*)]
\item $\beta(0) = 0$, $\beta'(0) = 1$, $\beta^{(\text{2l})}(0) = 0$ for $l \in \N$ and
\item $\beta|_{(0,\infty)} > 0$,
\end{enumerate}
the metric $g^{\beta}$ uniquely extends to a smooth rotationally symmetric metric on $\R^q$.

The pull-back of the curvature operator of $g^{\beta}$ is given by (cf. \cite[p.121]{petersen-2016})
\begin{equation}\label{eq:curvature-operator-beta}
  R^{\beta, q}(\theta, r) = \frac{1 - \beta'(r)^2}{\beta(r)^2} R_{\mathbb E \times \mathbb S^{q-1}} - \frac{\beta''(r)}{\beta(r)} L_q,
\end{equation}
where $L_q(e_i \wedge e_j) = \begin{cases} e_i \wedge e_j & \text{ if } i = 1, j > 1 \\ 0 & \text{ otherwise}\end{cases}$.

If the dimension $n$ is fixed, we will abuse notation denoting by $R^{\beta, q}$ and $L_q$ the curvature operators of $\R^n$ which are zero on combinations of the first $n-q$ co-ordinate directions.

\begin{dfn}\label{dfn:deformability}
  Let  $C \subset \mathcal C_{\mathrm B}(\E^n)$ be a codimension $c$ surgery stable curvature condition.
  It is called \emph{deformable}, if
  \begin{enumerate}[topsep=0cm]
  \item $0 \notin C$,
  \item it satisfies an inner ray condition with respect to $L_q$ for all $q \geq c$, i.e.\ $R \in C$ implies $R + \lambda L_q \in C$ for all $\lambda \geq 0$,
  \item and $\mu R_{\mathbb E^{n-q+1} \times \mathbb S^{q-1}} \in C$ for all $\mu > 0$ and $q \geq c$.
  \end{enumerate}
\end{dfn}

\begin{bsp}\label{bsp:deformable-curvature-conditions}
  \begin{enumerate}[label=(\roman*),leftmargin=1cm]
  \item By \cref{bem:curvature-condition-open-convex} \ref{item:bem:curvature-condition-open-convex-i}, a curvature condition $C$ given by an open convex cone with $R_{\mathbb E^{n-c+1} \times \mathbb S^{c-1}} \in C$ and $L_q \in C$ for all $n \geq q \geq c$ is deformable.
    This is the case for positive scalar curvature.
  \item Positive $p$-curvature for $0 \leq p < n-2$ is deformable (even though it does not satisfy an inner cone condition with respect to $L_q$ for $1 \leq p \leq n-2$).
  \item The condition $k$-positive Ricci curvature for $2 \leq k \leq n$ is deformable.
  \end{enumerate}
\end{bsp}

\begin{bem}
  Let $C \subset \mathcal C_{\mathrm B}(\E^n)$ be a deformable, codimension $c$ surgery stable curvature condition.
  Then $C$ satisfies an inner cone condition with respect to every $R_{\mathbb E^{n-q+1} \times \mathbb S^{q-1}} + \lambda L_q$ for all $q \geq c$ and $\lambda \geq 0$.
\end{bem}

\begin{prop}
  Suppose $C$ is a deformable, codimension $c$ surgery stable curvature condition.
  Then $R^{\beta, q} \in C$ for all $q \geq c$, if in $(0,\infty)$
  \begin{equation}\label{eq:rot-sym-satisf-C}
    1 - \beta'^2 > 0,
    \quad
    \beta'' \leq 0
    \quad \text{ and } \quad
    \beta'''(0) < 0.
  \end{equation}
\end{prop}

\begin{proof}
  From \eqref{eq:curvature-operator-beta} we conclude that, if
  \begin{align*}
    \mu \coloneqq \frac{1-\beta'^2}{\beta^2} > 0
    \quad \text{ and } \quad
    \lambda \coloneqq -\frac{\beta''}{\beta} \geq 0
  \end{align*}
  then by definition of deformability $R^{\beta, q} = \mu R_{\mathbb E^{n-q} \times \mathbb S^{n-1}} + \lambda L_q \in C$.
  By L'Hôpital, we have $\lim_{r \to 0}\mu = \lim_{r \to 0} \lambda = - \beta'''(0) > 0$ and thus $\lim_{r \to 0} R^{\beta, q} \in C$.
\end{proof}

\begin{dfn}
  One function of particular interest satisfying these assumptions is given by a smoothing of
  \begin{align*}
    \beta_{\delta} \colon [0, \infty) \to [0, \infty),
    \quad
    r \mapsto
    \begin{cases}
      \delta \sin(\nicefrac{r}{\delta}) & \text{ for } r \leq \frac{\delta \pi}{2}\\
      \delta & \text{ otherwise }
    \end{cases}
  \end{align*}
  where
  $\delta > 0$.
  The result, which we will also denote by $\beta_{\delta}$ is called a \emph{torpedo function of radius $\delta$}.
\end{dfn}

\begin{bem}
  Observing that $R_{\mathbb E^{n-q} \times \mathbb S^q} = R_{\mathbb E^{n-q} \times \mathbb E \times \mathbb S^{q-1}} + L_q$, we see that
  \begin{align*}
    R^{\beta_{\delta}, q}(\theta, r) =
    \begin{cases}
      R_{\mathbb E^{n-q} \times \mathbb S^q} & \text{ if } r = 0 \\
      R_{\mathbb E^{n-q} \times \mathbb E \times \mathbb S^{q-1}} & \text{ if } r \geq \frac{\delta \pi}{2}
    \end{cases}
  \end{align*}
  i.e.\ the metric induced by this function agrees on the last $q$ coordinates with a round metric near zero and has a cylindrical shape for large radii.
\end{bem}

\begin{dfn}\label{dfn:torpedo-metric}
  For every $q \geq 3$, the torpedo function $\beta_{\delta}$ gives rise to a metric $g^{\beta_{\delta}}$ on $\R^q$ and consequently on $D^q(r)$ (where we assume $r \geq \frac{\delta \pi}{2}$), which we will refer to as a \emph{torpedo metric} and denote by $g^{\delta}_{\torp}$.
\end{dfn}

\begin{bem}
  If $C \subset \mathcal C_{\mathrm B}(\E^n)$ is a codimension $c$ surgery stable curvature condition, and $n \geq q \geq c$, then $(\R^{n-q} \times \R^q, g_{\operatorname{eucl}} + g_{\torp}^{\delta})$ satisfies $C$.
\end{bem}

\begin{optional}
  \begin{dfn*}
    Let $\beta \colon [0,\delta] \to \R$ be a torpedo function.
    The arc-length parame\-trized curve $\gamma^{\beta}_{\torp}$ in $\R^2$ with $\pr_2 \circ \gamma^{\beta}_{\torp} \equiv \beta$ is called \emph{torpedo curve}.
    It induces an embedding, which is given in spherical co-ordinates by
    \begin{align*}
      i^{\beta}_{\torp} \colon D^q(\delta) \to \R \times \R^q,
      \quad
      (r, \theta) \mapsto (\alpha(r), \beta(r) \theta).
    \end{align*}
  \end{dfn*}

  \begin{prop*}
    We have $g^{\beta}_{\torp} = (i^{\beta}_{\torp})^{*} g_{\operatorname{eucl}}$.
  \end{prop*}
\end{optional}

\subsection{Connection metrics and riemannian submersions}\label{sec:connection-metrics}

Using the associated bundle construction, it is well-known that we can endow the total space of a vector bundle with a metric rotationally symmetric around the zero section.

\begin{prop}[{cf. \cite[Proposition 2.7.1, p.97]{gromoll-walschap-2009}}]\label{prop:connection-metrics}
  Let $(N,g_N)$ be a closed riemannian manifold and let $\pi \colon E \to N$ be a riemannian vector bundle of rank $q$ equipped with a metric connection $\omega$.
  Let $g_{\rot} = \dop r^2 + \beta^2(r)g_{\mathbb S^{q-1}}$ be a complete rotationally symmetric metric on $\R^q$.
  Then there exists a unique complete riemannian metric $h^{\nabla}$ on $E$ such that $\pi \colon (E,h^{\nabla}) \to (N,g_N)$ is a riemannian submersion with totally geodesic fibres isometric to $(\R^q, g_{\rot})$ and with horizontal distribution determined by $\omega$.
\end{prop}

\begin{dfn}
  We refer to $h^{\nabla}$ as a \emph{connection metric} on $E$ and write
  \begin{align*}
    h^\nabla = g_N \oplus_\omega g_{\rot}.
  \end{align*}
\end{dfn}

Recall that for a riemannian submersion $\pi \colon (E^n,g_E) \to (N^{n-q}, g_N)$ we can deform the metric $g_E$ by \emph{shrinking the fibre}.
We obtain a continuous path of metrics $\{g^t_E\}_{t \in (0,1]} \subset \mathcal R(E)$ given by
\begin{align*}
  g^t_E(X, Y) \coloneqq t^2 g_{F(\pi(p))}(X^{\mathcal V}, Y^{\mathcal V}) + \pi^{*}g_N(X,Y) \quad \text{ for } X,Y \in \Tangent_pE,
\end{align*}
where $X^{\mathcal V}, Y^{\mathcal V} \in \Tangent_{\pi(p)}F(\pi(p))$ are the orthogonal projections onto the tangent space at the fibre.
Clearly, we have $g^1_E = g_E$.

Because riemannian submersions from a complete riemannian manifold into any other riemannian manifold are fibre bundles (cf.\ \cite[Theorem 9.42 p.245]{besse-1987}), we can talk about the fibre of a riemannian submersion.
As a minor adaptation of \cite[Theorem 3.1]{hoelzel-2016}, we obtain the following result.

\begin{prop}\label{prop:hoelzel-scaling-submersions}
  Let $(E^n,g_E)$ is a complete riemannian manifold and let $C \subset \mathcal C_{\mathrm B}(\E^n)$ be a curvature condition.
  Further, let $\pi \colon (E^n,g_E) \to (N^{n-q}, g_N)$ be a riemannian submersion with fibres $\R^q$ into a closed manifold $N$.
  If $C$ satisfies an inner cone condition with respect to each curvature operator corresponding to
  \begin{align*}
    R_{b, p} := R_{( \R^{n-q} \times F(b), g_{\mathbb E^k} + g_E|_{F(b)})}(p),
  \end{align*}
  for all $b \in N$, $p \in F(b) \coloneqq \pi^{-1}(b) \cong \R^q$ and $R_{b, p}$ is constant on $\{ p \in E \mid d_{g_E}(p,N) > R \}$ for some $R > 0$, then there exists a $t_{*} > 0$ such that $g^t_E \in \mathcal R_C(E)$ for all $t \in (0,t_{*})$.
\end{prop}

\begin{optional}
  \begin{proof}
    We let $D := \{ p \in E \mid d_{g_E}(p,N) \leq R \}$ and apply the same reasoning as in \cite[Theorem 3.1]{hoelzel-2016}:
    Pick $U \subset D \subset E$ open with a $g_E$-o.n.\ frame $H(q) \oplus V(q)$ split in horizontal and vertical vectors.
    Set $V^t := \frac{1}{t}V$; then $H(q) \oplus V^{t}(q)$ is an o.n.\ frame for $g_E^t$.
    Write $R^t_E = R^t_F + E^t$ on $U$ (w.r.t.\ the frames, $R_F^t$ is the operator corresponding to $(e_1, \ldots, e_{n-q}) \oplus V$).

    Step 1.
    We know that for some $O(n)$-invariant cone $\tilde C$ and some compact $K$, we have $\tilde C\setminus K \subset C$.
    Moreover, for every $p \in U$ there exists a $\varepsilon_p > 0$ such that $B_{\varepsilon_p}(R_F^1(p)) \subset \tilde C$ (by the inner cone condition).

    Step 2.
    Now fix $p$ and $\varepsilon = \varepsilon_p$.
    Since $\tilde C$ is a cone, we have for all $t > 0$:
    \begin{align*}
      B_{\nicefrac{\varepsilon}{t^{2}}}(R^t_F(p)) \subset \tilde C
    \end{align*}
    and in particular $\exists t_p > 0: B_{\nicefrac{\varepsilon}{t_p^{2}}}(R^{t_p}_F(p)) \subset \tilde C\setminus K \subset C$.

    Step 3.
    Shrink $U$ (if necessary) to find $\delta > 0$ such that for all $q \in U$ and $t \in (0,t_p)$:
    \begin{align*}
      B_{\nicefrac{\delta}{t^{2}}}(R^t_F(q)) \subset C.
    \end{align*}
    Thus it is enough to find a $t_{*} \in (0, t_p]$ with $\|E^t\| \leq \nicefrac{\delta}{t^2}$ for all $q \in U, t \in (0,t_{*})$.

    Step 4.
    Check explicitly (\cite[Lemma 3.3.]{hoelzel-2016}) that $\|E^t\| \leq \nicefrac{C}{t}$ for some constant $C > 0$ and conclude $\exists t_{*}$ as above.
    Hence we have shown that $R^t_E \in C$ for $t \in (0,t_{*})$.

    Step 5.
    Also apply steps 1--4 to $U^{\infty} := E \setminus D$ where in step 3, we do not need to shrink, as $R_F^t$ does not depend on $q \in U$.

    Step 6.
    Find a finite cover of $E$ by open sets $U$ together with $U^{\infty}$ and conclude that we can choose $t_{*}$ such that $R^t_E$ is in contained in $C$ everywhere.
  \end{proof}
\end{optional}

Now consider the case that $g_E$ is a connection metric obtained from a rotationally symmetric metric $g_{\rot} = \dop r^2 + \beta^2(r) g_{\mathbb S^{q-1}}$ on $\R^q$.
Then shrinking the fibre amounts to shrinking the warping function, since at $b \in N$
\begin{align*}
  g_E^t|_{F(b)} = \dop r^2 + t^2 \beta^2(\textstyle\frac{r}{t}) g_{\mathbb S^{q-1}}.
\end{align*}
If $\beta$ is a torpedo function, the new warping function $r \mapsto t\beta(\frac{r}{t})$ can easily seen to be a torpedo function again.

\begin{kor}\label{kor:scaling-warped-metrics}
  Let $C \subset \mathcal C_{\mathrm B}(\E^n)$ be a deformable, codimension $c$ surgery stable curvature condition.
  Let $(N^{n-q},g_N)$ be a riemannian manifold, let $\pi \colon E \to N$ be a riemannian vector bundle of rank $q$ for $q \geq c$ equipped with a metric connection $\omega$ and let $g_{\rot} = \dop r^2 + \beta^2 g_{\mathbb S^{q-1}}$ be a rotationally symmetric metric on $\R^q$, which satisfies the conditions \cref{eq:rot-sym-satisf-C} and $\beta|_{[R, \infty)}$ is constant for some $R > 0$.
  Then there exists a $t_{*} > 0$ such that $(h^{\nabla})^t = g_N \oplus_{\omega} (\dop r^2 + t^2\beta^2(\frac{r}{t}) g_{\mathbb S^{q-1}})$ satisfies $C$ for all $t \in (0,t_{*})$.
\end{kor}

\begin{proof}
  Since $g_{\rot}$ satisfies the conditions \cref{eq:rot-sym-satisf-C}, we conclude that the curvature operators $R^{\beta,q}$ are contained in $C$ and by deformability, $C$ satisfies an inner cone condition with respect to them.
  These are precisely the curvature operators of the metric $g_{\mathbb E^{n-q}} + g_E|_{F(b)}$, where $F(b)$ denotes the fibre at $b \in N$.
  Applying \cref{prop:hoelzel-scaling-submersions} finishes the proof, because shrinking the metric is a deformation through connection metrics, as mentioned above.
\end{proof}

In particular, we can construct a connection metric on the total space of a vector bundle of suitable rank, which satisfies $C$, e.g.\ by considering the connection metric obtained from $g_{\torp}^{\delta}$ for $\delta$ small enough.

\subsection{Rotational symmetry around a submanifold}

Applying this to the normal bundle of a submanifold $N$ in $M$ we can produce a metric satisfying $C$ in a tubular neighbourhood of $N$ from a rotationally symmetric one.
If we start with an arbitrary metric in $M$, however, clearly it will not necessarily be a connection metric with respect to a rotationally symmetric metric on the normal bundle.
Nevertheless, it was observed by Gromov and Lawson that every metric does actually look almost rotationally symmetric close to $N$.

Throughout this section, let $M^n$ be a smooth manifold of dimension $n$, $i \colon N^k \hookrightarrow M^n$ be a closed submanifold with normal bundle $\pi \colon \nu N \to N$ equipped with a bundle metric $h^{\nu N}$.
Moreover, let $\phi \colon \nu N \to N$ be a tubular map, i.e.\ $\phi$ is an embedding with $\phi|_0 \equiv i \circ \pi$ and $\phi \circ s_0 \simeq i$ (where $s_0 \colon N \to \nu N$ is the zero section).

\begin{dfn}
  A riemannian metric $g$ on $M$ is called \emph{adjusted to the tubular map $\phi$ on the $r$-tube}, if $[0,r] \to M,\ r \mapsto \phi(r \nu_p)$ is a unit speed geodesic w.r.t.\ $g$, where $\nu_p \in \nu M$ with $\|\nu_p\|_{h^{\nu N}} = 1$ and $\|\blank\|_{h^{\nu N}}$ denotes the norm given by the riemannian vector bundle structure on $\nu N$.
\end{dfn}

From now on fix a metric connection $\omega$ on $\nu N$ and we denote by $\nu^{\leq r} N$ the radius $r$ disc bundle w.r.t.\ $h^{\nu N}$.

We can adjust an entire family of metrics to the tubular map $\phi$.

\begin{prop}[{\cite[adapted from Proposition 3.4]{ebert-frenck-2018}}]\label{prop:adjusting-to-phi}
  Let $\{g_{\xi}\}_{\xi \in D^l}$ be a continuous family of metrics such that $g_{\xi}$ for $\xi \in S^{l-1}$ is adjusted to $\phi$ on the $r$-tube.
  Then there exists an $r_0 \in (0, r]$ and a continuous map $F \colon [0,1] \times D^l \to \Diff(M)$ such that
  \begin{enumerate}[label=(\roman*)]
  \item $F|_{\{0\} \times D^l \cup [0,1] \times S^{l-1}} \equiv \id_M$,
  \item $F(t,x)|_N \equiv \id_N$ for all $(t,x) \in [0,1] \times D^l$,
  \item $(F(1,x))^{*}g_{\xi}$ is adjusted to $\phi$ on the $r_0$-tube.
  \end{enumerate}
\end{prop}

Let $C \subset \mathcal C_{\mathrm B}(\E^n)$ be a deformable, codimension $c$ surgery stable curvature condition and fix an arbitrary riemannian metric $g_N$ on the submanifold $N$.
By \cref{kor:scaling-warped-metrics}, we know that there exists a connection metric $h^{\torp} := g_N \oplus_{\omega} g_{\torp}^{\delta}$, which satisfies $C$ and which we will also fix from now on.

\begin{dfn}\label{dfn:metrics-in-standard-form}
  Let $R > 0$ be a fixed radius.
  We call a riemannian metric $g \in \mathcal R_C(M)$ \emph{rotationally symmetric around $N$}, if
  \begin{enumerate}
  \item\label{dfn:metrics-in-standard-form:1} $g$ is adjusted to $\phi$ on the $R$-tube,
  \item\label{dfn:metrics-in-standard-form:2} $\phi|_{\nu^{\leq R}N}^{*}g = (g_N \oplus_{\omega} g_{\rot})|_{\nu^{\leq R}N}$ for some rotationally symmetric metric $g_{\rot}$ on $\R^{n-k}$ for which $(\R^n, g_{\operatorname{eucl}} + g_{\rot})$ satisfies $C$.
  \end{enumerate}
  Now define the \emph{space of rotationally symmetric metrics around $N$} as
  \begin{align*}
    \mathcal R^{\rot}_C(M) & \coloneqq \mathcal R^{\rot}_C(M; \phi, g_N, h^{\nu N}, \omega, R)\\
                                         & \coloneqq \{ g \in \mathcal R_C(M) \text{ is rotationally symmetric around } N \}.
  \end{align*}
  Moreover, we denote by
  \begin{align*}
    \mathcal R^{\torp}_C(M) & \coloneqq \mathcal R^{\torp}_C(M; \phi, h^{\nu N}, h^{\torp}, R)\\
                            & \coloneqq \{ g \in \mathcal R_C(M) \mid \phi|_{\nu^{\leq R}N}^{*}g = h^{\torp}|_{\nu^{\leq R}N} \}.
  \end{align*}
  the \emph{space of riemannian metrics satisfying $C$, which are standard near $N$}.
\end{dfn}

\begin{bem}
  \begin{enumerate}[label=(\roman*)]
  \item The adjustment in \cref{dfn:metrics-in-standard-form:1} simply amounts to $g_{rot}$ in \cref{dfn:metrics-in-standard-form:2} being of the form $\dop r^2 + \beta g_{\mathbb S^{n-1}}$ in terms of Fermi coordinates around $N$.
  \item Note that metrics in $\mathcal R^{\torp}_C(M)$ are adjusted to $\phi$ on the $R$-tube by definition and thus are rotationally symmetric.
    In particular, we have an inclusion $\mathcal R^{\torp}_C(M) \hookrightarrow \mathcal R^{\rot}_C(M)$.
  \end{enumerate}
\end{bem}

\section{Main Results}\label{section:main-results}

\subsection{Main technical result}

We will now state the main technical result of this work, which shows that the space of all metrics satisfying curvature conditions of a certain type is weakly homotopy equivalent to the space of metrics, which take a particular prescribed form around an initially fixed submanifold.
The proof of this theorem is heavily built on the techniques and terminology developed in \cite{chernysh-2004,walsh-2011} and essentially follows Chernysh's presentation implementing some improvements introduced by \cite{ebert-frenck-2018}.

\begin{satz}\label{satz:c}
  Let $C \subset \mathcal C_{\mathrm B}(\E^n)$ be a deformable, codimension $c$ surgery stable curvature condition.
  \begin{enumerate}
  \item Let $M^n$ be a smooth manifold of dimension $n$ and let $N^k$ be a compact, $k$-dimensional submanifold in $M$ with $n - k \geq c$.
  \item Let $\phi \colon \nu N \to M$ be tubular map and let $h^{\nu N}$ be a bundle metric on the normal bundle $\nu N$ equipped with a metric connection $\omega$.
  \item Further let $g_N$ be an arbitrary riemannian metric on $N$ and let $h^{\torp} = g_N \oplus_{\omega} g_{\torp}^{\delta}$ be a connection metric on $\nu N$, which satisfies $C$ (as constructed in \cref{sec:connection-metrics}), obtained from $g_N$, $h^{\nu N}, \omega$ and a torpedo metric $g_{\torp}^{\delta}$ on $\R^{n-k}$.
  \end{enumerate}
  Then the inclusion of metrics, which are standard near $N$
  \begin{align*}
    i \colon \mathcal R^{\torp}_C(M) \hookrightarrow \mathcal R_C(M)
  \end{align*}
  is a weak homotopy equivalence.
\end{satz}

The proof of \cref{satz:c} is based on the following two propositions whose proofs will occupy the rest of this paper.

\begin{prop}\label{prop:deformation-of-families}
  Let $g \colon S \to \mathcal R_C(M),\ \xi \mapsto g_{\xi}$ be a family of metrics.
  There exists a continuous map
  \begin{align*}
    \Pi \colon [0,1] \times S \to \mathcal R_C(M)
  \end{align*}
  with the following properties:
  \begin{enumerate}[label=(\roman*)]
  \item \ $\Pi(0, \blank) \equiv g$,
  \item \label{item:prop:def:ii} \ $\Pi(\{1\} \times S) \subset \mathcal R^{\rot}_C(M)$,
  \item \label{item:prop:def:iii} \ If $g_{\xi} \in \mathcal R^{\rot}_C(M)$, then $\Pi(t, \xi) \subset \mathcal R^{\rot}_C(M)$ for all $t \in [0,1]$.
  \end{enumerate}
\end{prop}

The proof will be carried out in \cref{section:proof-satz-b}.

\begin{prop}\label{prop:deformation-warped-to-torpedo}
  The space $\mathcal R^{\torp}_C(M)$ is a weak deformation retract of $\mathcal R^{\rot}_C(M)$.
\end{prop}

The proof will be carried out in \cref{section:rotationally-sym-metrics}.

\begin{proof}[Proof of \cref{satz:c}]
  By \cref{prop:deformation-warped-to-torpedo}, it is enough to show that the inclusion $i \colon \mathcal R_C^{\rot}(M) \hookrightarrow \mathcal R_C(M)$ is a weak homotopy equivalence.

  It is well-known that $i$ is a weak homotopy equivalence, if for every $n \in \N_0$ and map $g \colon \{0\} \times D^n \to \mathcal R_C(M)$ with $g(\{0\} \times S^{n-1}) \subset \mathcal R_C^{\rot}(M)$ there exists a homotopy $\overline g \colon [0,1] \times D^n \to \mathcal R_C(M)$ such that the following diagram commutes
  \begin{equation*}
    \begin{tikzcd}[row sep=small]
      \{0\} \times S^{n-1} \ar{rr}{g} \ar[hook]{dr} \ar[hook]{ddd} & & \mathcal R_C^{\rot}(M) \ar[hook]{ddd}{i} \\
      & {[0,1] \times S^{n-1} \cup \{1\} \times D^n \ar[hook]{d} \ar[dashed]{ru}{\overline g}} & \\
      & {[0,1] \times D^n \ar[dashed]{dr}{\overline g}} & \\
      \{0\} \times D^n \ar{rr}{g} \ar[hook]{ru} & & \mathcal R_C(M)
    \end{tikzcd}
  \end{equation*}
  Applying \cref{prop:deformation-of-families} to $g$, we let $\overline g \coloneqq \Pi$.
  By \cref{item:prop:def:ii}, $\overline g(\{1\} \times D^n) \subset \mathcal R_C^{\rot}(M)$ and since $\overline g(\{0\} \times S^{n-1}) \subset \mathcal R_C^{\rot}(M)$, we conclude with \cref{item:prop:def:iii} that $\overline g([0,1] \times S^{n-1}) \subset \mathcal R_C^{\rot}(M)$.
\end{proof}

\subsection{Applications}

The following is a theorem of Chernysh \cite[Theorem 1.2]{chernysh-2004} and Walsh \cite[Theorem 4.1]{walsh-2013} for the case of $C = \psc$.
It is the technical version of and clearly implies \cref{msatz:homotopy-invariance}.

\begin{kor}\label{satz:a}
  Let $C \subset \mathcal C_{\mathrm B}(\E^n)$ be a deformable, codimension $c$ surgery stable curvature condition and let $M$ be a closed smooth manifold of dimension $n$.
  If $\chi(M,\phi)$ is obtained from $M$ by surgery along $\phi|_{S^k \times \{0\}}$ for $\phi \in \Emb(S^k \times D^{n-k}, M)$ for $k \geq c - 1$ and $n - k \geq c$, then
  \begin{align*}
    \mathcal R_C(M) \simeq \mathcal R_C(\chi(M, \phi)).
  \end{align*}
\end{kor}

It results easily from \cref{satz:c} applied to a surgery and its reversal, combined with the observation that the spaces $\mathcal R_C(M)$ are open subsets of a Fréchet manifold.
Thus, by work of Palais \cite[Theorem 14]{palais-1966}, they are dominated by CW-complexes and we obtain the homotopy equivalence in \cref{satz:a} from Whitehead's theorem.

\begin{optional}
  \begin{proof}[Proof of \cref{satz:a}]
    If $\chi(M,\phi)$ is obtained from $M$ by surgery of codimension at least $c$ for $\phi \colon S^k \times D^{n-k} \to M$, we know that there exists an embedding $\psi \colon D^{k+1} \times S^{n-k-1} \to \chi(M, \phi)$ with $\psi|_{S^k \times S^{n-k-1}} \equiv \phi_{S^k \times S^{n-k-1}}$ such that $M \setminus \Im \phi$ is diffeomorphic to $\chi(M,\phi) \setminus \Im \psi$.

    Now observe that the spaces of metrics satisfying $C$ on $M$ and $\chi(M,\phi)$, which are fixed within $\mathring{\Im \phi}$ or $\mathring{\Im \psi}$ resp., and have a product shape near $\phi(S^k \times S^{n-k-1})$ are homeomorphic.
    In particular, if we choose any metrics $g_{S^k}$, $g_{S^{n-k-1}}$ on $S^k$ and $S^{n-k-1}$ and torpedo metrics $g_{\torp}^1$, $g_{\torp}^2$ on $D^{n-k}$ and $D^{k+1}$, the spaces of metrics $\mathcal R^{\torp}_C(M; \phi, g_{S^k}, g_{\torp}^1)$ and $\mathcal R^{\torp}_C(\chi(M, \phi); \psi, g_{S^{n-k-1}}, g_{\torp}^2)$ are homeomorphic.

    Thus, invoking \cref{satz:c} twice shows that $\mathcal R_C(M)$ is indeed homotopy equivalent to $\mathcal R_C(\chi(M,\phi))$.
  \end{proof}
\end{optional}

\begin{bem}
  The conditions $k \geq c-1$ and $n-k \geq c$ are needed to reverse the surgery and again have a sufficiently high codimension.
  Combined we see that $n-c \geq k \geq c-1$ and conclude that the surgery stability needs to be roughly ``below the middle dimension'' if we want to apply \cref{satz:a}.
\end{bem}

We note that we can slightly generalize the definition of a surgery as follows.
Let $W_1$ and $W_2$ be manifolds with boundary with diffeomorphic and connected boundaries $\partial W_1 \cong \partial W_2$.
If $\phi \colon W_1 \hookrightarrow M$ is an embedding, we can consider
\begin{align*}
  \tilde \chi(M; \phi) \coloneqq M \setminus \operatorname{Int}(\phi(W_1)) \cup_{\partial W_1} W_2.
\end{align*}

We can apply this procedure to prove \cref{msatz:hpn}, which follows from the following corollary since both positive scalar curvature and positive $1$-curvature are deformable and  admit codimension $4$ surgeries (cf.\ \cref{bsp:ss-cc-psc,bsp:ss-cc-p-curv,bsp:deformable-curvature-conditions}).

Denote by $g^{\delta}$ for $\delta \in (0,1]$ the metric on $S^{4k+3}$ obtained via shrinking the fibres of the riemannian submersion given by the Hopf fibration $S^3 \hookrightarrow \mathbb S^{4k+3} \to \mathbb HP^k$.

\begin{kor}\label{kor:space-of-metrics-S-HP}
  Let $k \geq 1$, let $C \subset \mathcal C_{\mathrm B}(\E^{4(k+1)})$ be a deformable curvature condition, which admits codimension $4$ surgeries and such that $g^{\delta} + \dop t^2$ satisfies $C$ for all $\delta \in (0,1]$.
  Then spaces of metrics satisfying $C$ on $S^{4(k+1)}$ and $\mathbb HP^{k+1}$ have the same homotopy type.
\end{kor}

A key ingredient in the proof of this statement is Gajer's Lemma, which was originally stated for positive scalar curvature, but holds in greater generality.

\begin{lem}[{\cite[p.184]{gajer-1987}}]\label{lemma:gajer}
  Let $C \subset \mathcal C_{\mathrm B}(\mathbb E^n)$ be a curvature condition, let $N^{n-1}$ be a closed manifold and let $\{g_t\}_{t \in [0,1]} \subset \mathcal R(N)$ be a smooth path of riemannian metrics.
  If for every $t \in [0,1]$ the riemannian product metrics $(N, g_t) \times \mathbb E$ satisfy $C$, then there exists a $0 < \Lambda \leq 1$ such that for every smooth function $f \colon \R \to [0,1]$ with $|f'|, |f''| \leq \Lambda$ the metric $g_{f(t)} + \dop t^2$ on $N \times \R$ satisfies $C$.
\end{lem}

\begin{proof}
  The proof follows directly from calculations analogous to \cite[p.185]{gajer-1987}, which yield
  \begin{align*}
    R_{(N \times \R, g_{f(t)} + \dop t^2)} & = R_{(N, g_{f(t_0)}) \times \mathbb E} + O(|f'|) E_1 + O(|f'|^2) E_2 + O(|f''|) E_3,
  \end{align*}
  where $E_1, E_2, E_3$ only depend on the family of metrics $\{g_t\}$ and we note that $\tr E_1 = 0$.
  Since $C$ is an open subset in $\mathcal C_{\mathrm B}(\mathbb E^n)$, we find a suitable choice for $\Lambda$.
\end{proof}

\begin{proof}[Proof of \cref{kor:space-of-metrics-S-HP}]
  Consider the inclusion $i \colon \mathbb HP^{k} \hookrightarrow \mathbb HP^{k+1}$.
  It is well-known that the normal bundle $\nu \mathbb HP^k$ can be written as an associated bundle $S^{4k+3} \times_{S^3} \R^4$ to the Hopf fibration, i.e.\ to the $S^3$-principal bundle $S^{4k+3} \to \mathbb HP^k$.

  Now choose a bundle metric $h^{\nu \mathbb HP^k}$ and $\overline r > 0$ small enough such that the inclusion of the disc bundle $\phi \colon \nu^{\leq \overline r} \mathbb HP^k \hookrightarrow \mathbb HP^{k+1}$ is an embedding.
  Then we have $\partial(\nu^{\leq \overline r} \mathbb HP^k) = \nu^{\overline r} \mathbb HP^k \cong S^{4k+3}$ and it is well-known that the complement of the disc bundle is diffeomorphic to $D^{4k}$.
  Thus we conclude that
  \begin{align*}
    \big(\mathbb HP^{k+1} \setminus (\nu^{\leq \overline r}\mathbb HP^k)\big) \cup_{S^{4k+3}} D^{4(k+1)} \cong S^{4(k+1)}
  \end{align*}
  In the same vain, we obtain
  \begin{align*}
    \big(S^{4(k+1)} \setminus D^{4(k+1)}\big) \cup _{S^{4k+3}} \big(S^{4k+3} \times_{S^3} D^4\big) \cong \mathbb HP^{k+1}.
  \end{align*}
  In both cases we remove a tubular neighbourhood of a codimension $4$ submanifold and thus we are in a position to apply \cref{satz:c} to both situations.
  It remains to check that the spaces of metrics which are standard near the submanifolds $\mathcal R^{\torp}_C(\mathbb HP^{k+1})$ and $\mathcal R^{\torp}_C(S^{4(k+1)})$ are weakly homotopy equivalent for suitable choice of torpedo- and connection metrics.
  We can strengthen the resulting weak homotopy equivalence between $\mathcal R_C(\mathbb HP^{k+1})$ and $\mathcal R_C(S^{4(k+1)})$ to a homotopy equivalence using Palais and Whitehead, as mentioned after \cref{satz:a}.

  Let $g^{\delta}_{\torp}$ be a torpedo metric.
  Then $S^3$ acts isometrically on $\mathbb S^{4k+3}$ from the right via the Hopf fibration and isometrically on $(\R^4, g^{\delta}_{\torp})$ from the left, by the torpedo metric's rotational symmetry.
  Now we obtain a riemannian metric on the quotient $S^{4k+3} \times_{S^3} \R^4$, induced from the product metric $g_{\mathbb S^{4k+3}} + g_{\mathbb E^4}$, such that the projection map to $\mathbb HP^k$ is a riemannian submersion.
  Because the Hopf fibration has totally geodesic fibres and the action of $S^3$ is by isometries, the riemannian submersion $S^{4k+3} \times_{S^3} \R^4 \to \mathbb HP^k$ has totally geodesic fibres (cf.\ \cite[p.98]{gromoll-walschap-2009}).
  But such metrics are already connection metrics \cite[Theorem 2.7.2, p.98]{gromoll-walschap-2009}.

  By \cref{kor:scaling-warped-metrics}, we can choose $\delta > 0$ small enough such that the connection metric on $S^{4k+3} \times_{S^3} \R^4$ satisfies $C$.
  Moreover, $(\R^4 \setminus D^4(\frac{\delta \pi}{2}), g^{\delta}_{\torp})$ is isometric to a riemannian product $\mathbb S^3(\delta) \times (\frac{\delta \pi}{2}, \infty)$ and thus we obtain isometries
  \begin{align*}
    \big(S^{4k+3} \times_{S^3} \R^4\big) \setminus \big(S^{4k+3} \times_{S^3} D^4(\nicefrac{\delta \pi}{2})\big)
    & = S^{4k+3} \times_{S^3} \big(\mathbb S^3(\delta) \times (\nicefrac{\delta \pi}{2}, \infty)\big)\\
    & = \big(S^{4k+3} \times_{S^3} \mathbb S^3(\delta)\big) \times (\nicefrac{\delta \pi}{2}, \infty).
  \end{align*}
  We denote the metric obtained from this on $S^{4k+3} \cong S^{4k+3} \times_{S^3} \mathbb S^3(\delta)$ by $g^{\delta}$.

  For any $0 < \delta \leq 1$ the metric $g^{\delta} + \dop t^2$ on $S^{4k+3} \times \R$ satisfies $C$ by assumption.

  Now choose a smooth function $s \colon [0,1] \to \R$, which is constant near $\{0,1\}$ and satisfies $s(0) = 1$, $s(1) = \delta$ and $s' \leq 0$.
  We obtain a smooth path of metrics $\{g^{s(t)}\}_{t \in [0,1]}$ on $S^{4k+3}$, which satisfies the assumptions of \cref{lemma:gajer}.

  If we choose the inclusion $i \colon \operatorname{pt} \hookrightarrow S^{4(k+1)}$ to a point in the sphere, the space $\mathcal R^{\torp}_C(S^{4(k+1)})$ is the space of metrics, which are fixed and of torpedo shape on one hemisphere, i.e. $\phi|_{\nu^{\leq R}\operatorname{pt}}^{*}g = h^{\torp}_{\nu^{\leq R}\operatorname{pt}}$.
  In particular, there exist $0 < R' \leq R$ and $b > 0$ such that for every $g \in \mathcal R^{\torp}_C(S^{4(k+1)})$ we have that $\phi|_{\nu^{R' \leq R}\operatorname{pt}}^{*}g$ is isometric to $(S^{4k + 3} \times [0, 2b], g_{\mathbb S^{4k + 3}} + \dop t^2)$.
  Possibly by passing to another torpedo metric $h^{\torp}$, we obtain $b$ large enough such that $g^{\frac{1}{b}s(t)} + \dop t^2$ is a metric satisfying $C$ on $S^{4k+3} \times [0,b]$, by \cref{lemma:gajer}.
  Denote by $\tilde R$ the space of metrics satisfying $C$ with $\phi|_{\nu^{\leq R'}\operatorname{pt}}^{*}g = h^{\torp}_{\nu^{\leq R'}\operatorname{pt}}$ and $\phi|_{\nu^{R'\leq R}\operatorname{pt}}^{*}g$ is isometric to
  \begin{align*}
    (S^{4k+3} \times [0, b], g_{\mathbb S^{4k + 3}} + \dop t^{2}) \cup (S^{4k+3} \times [0, b], g^{\frac{1}{b}s(t)} + \dop t^2).
  \end{align*}

  The spaces $\mathcal R^{\torp}_C(S^{4(k+1)})$ and $\tilde R$ are weakly homotopy equivalent, since we can use the family $\{g^{s(t)}\}_{t \in [0,1]}$ and Gajer's lemma to interchange the cylindrical pieces on $S^{4k+3} \times [0,2b]$.

  The spaces $\mathcal R^{\torp}_C(\mathbb HP^{k+1})$ and $\tilde R$ are homeomorphic.
\end{proof}

\begin{optional}
  \begin{bem*}
    Using the metric constructed in the proof, we can also construct a metric on $\mathbb HP^k \# \mathbb HP^k$, which satisfies a curvature condition admitting codimension $4$ surgeries.
    This procedure's reversal however is a (generalized) surgery of codimension $1$ and thus we cannot apply \cref{satz:c} in this case.
    Moreover, even admittance of codimension $4k$ surgeries in dimension $4k$ would allow us to construct a metric on the connected sum.
  \end{bem*}

  \begin{bem*}
    We would like to obtain a condition, such that the assumption on the metrics $g^{\delta}$, given by shrinking the fibres of the Hopf fibration $S^{3} \to \mathbb S^{4k+3} \xrightarrow{\pi} \mathbb HP^k$, follows.
    The formulas for the canonical variation (cf.\ e.g.\ \cite[Lemma 3.3]{hoelzel-2016}) can be applied.
    Since $\pi$ has totally geodesic fibres, the terms involving the $T$-tensor vanish.
    Falcitelli et.\ al\footnote{Falcitelli, Maria; Ianus, Stere; Pastore, Anna Maria. Riemannian submersions and related topics. World Scientific Publishing Co., Inc., River Edge, NJ, 2004), p.20} compute that for $X,Y \in \mathcal V^{h}(S^{4k+3})$
    \begin{align*}
      A_XY = - \sum_{i = 1}^{3} g_{\mathbb S^{4k+3}}(I_i X, Y)I_iN,
    \end{align*}
    where $I_1, I_2, I_3 \colon \Tangent \R^{4(k+1)} \to  \Tangent \R^{4(k+1)}$ are the three almost-complex structures given by the quaternionic structure in $\R^{4(k+1)}$ and $N$ is the unit normal ($I_iN$ are unit vector fields spanning the tangent bundle of the fibres).

    Let $g = g_{\mathbb S^{4k+3}}$, $X,Y \in \mathcal V^h(S^{4k+3})$, $U, V \in \mathcal V^v(S^{4k+3})$, then $A_XU \in \mathcal V^h(S^{4k+3})$ and $g(A_XU, Y) = -g(U,A_XY)$, thus
    \begin{align*}
      g(A_XU, A_YV)
      & = -g(U, A_XA_YV)
      = \sum_{i = 1}^3 g(I_iX,A_YV)g(U,I_iN)\\
      & = \sum_{i,j = 1}^3 g(I_iX, I_jY)g(I_iN, U)g(I_jN,V).
    \end{align*}
    Therefore, we have $g(A_XU, A_YV) = g(A_XV, A_YU)$, which causes the terms involving $A$ in the third row of \cite[Lemma 3.3]{hoelzel-2016}) to vanish.

    As in \cite[Proof of Theorem 3.1]{hoelzel-2016}, let $(v_1, v_2, v_3) \oplus (h_1, \ldots, h_{4k})$ be an orthonormal frame and decomposition w.r.t.\ $g_{\mathbb S^{4k+3}}$ and the submersion at some point $p \in S^{4k+3}$.
    Choose $v_i$ such that $v_i = (I_iN)|_p$.
    Then $(\frac{1}{t}v_1, \frac{1}{t}v_3, \frac{1}{t}v_3) \oplus (h_1, \ldots, h_{4k})$ is an orthonormal frame and decomposition w.r.t.\ $g_{\mathbb S^{4k+3}}^t$ and we have
    \begin{align*}
      \tilde R^t_{\mathbb S^{4k+3}}
      = \tilde R^t_{\mathbb S^3 \times \mathbb E^{4k}} + E^t
      = \frac{1}{t^2}\tilde R_{\mathbb S^3 \times \mathbb E^{4k}} + E^t
    \end{align*}
    with $i,j,m,l \in \{1, 2, 3\}$ and $r,s,u,v \in \{1, \ldots, 4k\}$
    \begin{align}
      E^t (e_i, e_j, e_m, e_l) & = 0 \nonumber\\
      E^t (e_i, e_j, e_m, e_{3+r}) & = t^{-1} R_{\mathbb S^{4k+3}}(v_i, v_j, v_m, h_r) - 0 = 0 \label{eq:ijmr}\\
      E^t (e_i, e_j, e_{3+r}, e_{3+s}) & = R_{\mathbb S^{4k+3}}(v_i, v_j, h_r, h_s) \label{eq:ijrs}\\
                               & \quad + \left(1 - t^2 \right) \underbrace{\left( g(A_{h_r}v_j, A_{h_s}v_i) - g(A_{h_r}v_i, A_{h_s}v_j) \right)}_{= 0 \text{ (see above)}} \nonumber \\
      & = 0 \nonumber\\
      E^t(e_{3+r}, e_i, e_{3+s}, e_j) & = R_{\mathbb S^{4k+3}}(h_r, v_i, h_s, v_j) + (1 - t^2) g(A_{h_r}v_j, A_{h_s}v_i) \label{eq:risj} \\
      E^t(e_{3+r}, e_{3+s}, e_{3+u}, e_i) & = t R_{\mathbb S^{4k+3}}(h_r, h_s, h_u, v_i) = 0 \label{eq:rsui} \\
      E^t(e_{3+r}, e_{3+s}, e_{3+u}, e_{3+v}) & = t^2 R_{\mathbb S^{4k+3}} (h_r, h_s, h_u, h_v)\\
      & \quad + (1 - t^2) (\pi^*R_{\mathbb HP^k})(h_r, h_s, h_u, h_v). \nonumber
    \end{align}
    Note that the contributions of $R_{\mathbb S^{4k+3}}$ in the \eqref{eq:ijmr},\eqref{eq:ijrs} and \eqref{eq:rsui} vanish, because $R_{\mathbb S^{4k+3}} = \id_{\bigwedge^2\mathbb E^{4k+3}}$, i.e. $\tilde R_{\mathbb S^{4k+3}}(e_i, e_j, e_m, e_l) = \langle e_i, e_l\rangle\langle e_j,e_m\rangle - \langle e_i,e_m\rangle\langle e_j,e_l\rangle = \delta_{il}\delta_{jm} - \delta_{im}\delta_{jl}$ for $i,j,m,l \in \{1, \ldots, 4k+3\}$.
    In \eqref{eq:risj}, we obtain from above
    \begin{align*}
      g(A_{h_r}v_j, A_{h_s}v_i) = \sum_{a,b = 1}^3g(I_ah_r, I_bh_s)g(v_a,v_j)g(v_b,v_i) = g(I_jh_r, I_ih_s).
    \end{align*}
    As we write
    \begin{align*}
      E^t(e_{3+r}, e_{3+s}, e_{3+u}, e_{3+v}) & = (R_{\mathbb S^{4k+3}}  - R_{\mathbb S^{4k+3}}  + t^2 R_{\mathbb S^{4k+3}}) (h_r, h_s, h_u, h_v),\\
      & \quad + (1 - t^2) (\pi^*R_{\mathbb HP^k})(h_r, h_s, h_u, h_v)
    \end{align*}
    we notice that we can express
    \begin{align*}
      E^t = \tilde R_{\mathbb S^{4k+3}} - \tilde R_{\mathbb S^3 \times \mathbb E^{4k}} + (1-t^2) \tilde E
    \end{align*}
    where (in the notation, as in the definition of $E$, all other values $=0$):
    \begin{align*}
      \tilde E(e_{3+r}, e_i, e_{3+s}, e_j) & = g(A_{h_r}v_j, A_{h_s}v_i) = g(I_jh_r, I_ih_s)\\
      \tilde E(e_{3+r}, e_{3+s}, e_{3+u}, e_{3+v}) & = (\pi^*R_{\mathbb HP^k})(h_r, h_s, h_u, h_v) - R_{\mathbb S^{4k+3}} (h_r, h_s, h_u, h_v).
    \end{align*}
    In total, we have
    \begin{align*}
      \tilde R^t_{\mathbb S^{4k+3}}
      = \tilde R_{\mathbb S^{4k+3}} + \frac{1 - t^2}{t^2} \tilde R_{\mathbb S^3 \times \mathbb E^{4k}} + (1-t^2)\tilde E \in \mathcal C_{\mathrm B}(\mathbb E^{4k+3}).
    \end{align*}
    Thus, we can guarantee that $R^t_{\mathbb S^{4k+3}} \in C$, if $\|\tilde E\| \leq \rho(R_{\mathbb S^3 \times \mathbb E^{4k}})$ for $C$.
  \end{bem*}
\end{optional}

\section[Proof of Theorem 3.4]{Proof of \cref{prop:deformation-of-families}}\label{section:proof-satz-b}

Before starting the actual proof, we will recall the classical graph deformation procedure introduced by Gromov and Lawson.
The construction of $\Pi$ then proceeds in two steps.
First, the graph deformation is applied to a family of riemannian metrics to split a tubular neighbourhood around a submanifold $N$ into three regions with particular properties.
Then one deforms the metric on these three regions to obtain a metric, which is rotationally symmetric around $N$.

\subsection{Preliminaries and Chernysh's trick}

Here we will recall a construction for a metric deformation and a few technical results of Hoelzel.

Before, let us recall elementary facts about curves in $\R^2$, which will be used to control the metric deformations.
We will deal with arc-length parametrized curves $\gamma \colon \R \to \R^2, s \mapsto (r(s), t(s))$, which satisfy a number of properties.
\begin{dfn}
  Let $\overline r > 0$.
  We denote by
  \begin{align*}
    \Gamma_{b}(\overline r) \coloneqq \{ \gamma \colon \R \to \R^2 \text{ arc-length parametrized curve satisfying \ref{item:gamma-near-0} -- \ref{item:gamma-b-last} below }\}
  \end{align*}
  \begin{enumerate}[label=(\roman*)]
  \item\label{item:gamma-near-0} $\gamma(0) = (\overline r, 0)$ and $t|_{(-\infty,0]} \equiv 0$,
  \item $t(s) \geq 0$ for all $s \in \R$,
  \item $\gamma$ intersects the $t$-axis $\{0\} \times \R$ precisely once following the arc of a circle (of possibly infinite radius) at $\gamma(b)$ and is symmetric about it,
  \item \label{item:gamma-b-last} $r$ is non-increasing, while $t$ is non-decreasing for $s \in (-\infty,b]$.
  \end{enumerate}
  \begin{align*}
    \tilde \Gamma_b(\overline r) \coloneqq \{ \gamma \in \Gamma_b(\overline r) \mid \gamma \text{ satisfies \ref{item:gamma-b-tilde} below } \}
  \end{align*}
  \begin{enumerate}[label=(\roman*),resume]
  \item \label{item:gamma-b-tilde} There exists a partition $0 = s_0 \leq s_1 \leq \cdots \leq s_6 = b$ such that
    \begin{align*}
      \kappa|_{[s_0,s_1] \cup [s_2,s_3] \cup [s_4,s_5]} \equiv 0
      \quad \text{ and } \quad
      r'|_{[s_3, s_4]} \equiv 0
    \end{align*}
    where $\kappa$ is the signed curvature function of $\gamma$.
  \end{enumerate}
  We endow each of the above sets with the subspace topology from $C^{\infty}(\R, \R^2)$ and denote $\Gamma(\overline r) \coloneqq \bigcup_{b > 0} \Gamma_b(\overline r)$ and $\tilde \Gamma(\overline r) \coloneqq \bigcup_{b > 0} \tilde \Gamma_b(\overline r)$.
\end{dfn}

\begin{figure}[h!]
  \centering
  \includegraphics[draft=false,width=.6\textwidth]{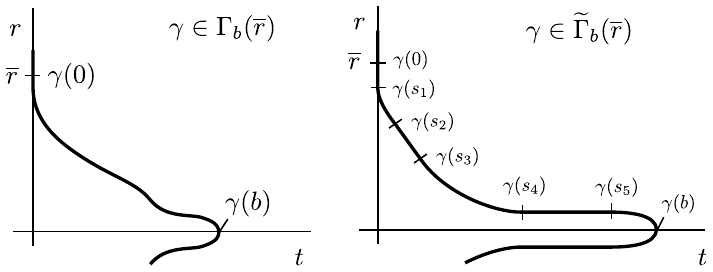}
  \caption{Examples for curves in $\Gamma_b(\overline r)$ and $\tilde\Gamma_b(\overline r)$.}
\end{figure}

\begin{optional}
  \begin{bem*}
    Clearly we have the following.
    \begin{enumerate}
    \item The map $L \colon \Gamma(\overline r) \to (0,\infty), \ \gamma \mapsto \gamma^{-1}(\{0\} \times \R)$, which determines the length of the curve segment between $(\overline r,0)$ and the curve's intersection with the $t$-axis, is continuous.
    \item The space $\Gamma(\overline r_1)$ is contained as a subspace in $\Gamma(\overline r_2)$ for $r_1 \leq r_2$.
    \end{enumerate}
  \end{bem*}
\end{optional}

\begin{prop}
  Any curve $\gamma$ in $\Gamma_b(\overline r)$ is uniquely and continuously determined on the interval $[0,b]$ by each of the following:
  \begin{enumerate}
  \item\label{item:determined-on-0-L} $\gamma|_{[0,b]}$,
  \item its angular function $\theta \colon \R \to [0,\frac{\pi}{2}]$,
  \item its signed curvature function $\kappa \colon \R \to \R$.
  \end{enumerate}
\end{prop}

\begin{proof}
  The claim \ref{item:determined-on-0-L} is immidiately clear.
  For the remaining note that the curve $\gamma$ is determined by its angular function $\theta \colon \R \to [0,\frac{\pi}{2}]$ as it gives rise to the following initial value problem
  \begin{align*}
    \begin{cases}
      \cos \theta = \left<\gamma', -\partial_r\right> = -r' ,
      & r(0) = \overline r\\
      (r')^2 + (t')^2 = 1
      & t(0) = 0\\
      t \geq 0
    \end{cases}.
  \end{align*}
  The angular velocity $\theta'$ is precisely the signed curvature of $\gamma$ and $\theta(0) = 0$ (by \ref{item:gamma-near-0}).
  Moreover, we conclude that by the theory of ordinary differential equations $\gamma$ continuously depends on $\theta$ or $\kappa$, respectively.
\end{proof}

Note that from this proof we can extract a description for $r(s)$ and $t(s)$ for $s \in [0,b]$ in terms of $\theta$ as follows
\begin{equation}
  r(s) = \overline r - \int_0^s \cos \theta(u) \dop u
  \text{ and }
  t(s) = \int_0^s \sin \theta(u) \dop u.
  \label{eq:r-t-in-terms-of-theta}
\end{equation}

If $\gamma$ is determined by $\theta$ or $\kappa$, we will write $\gamma(\theta)$ or $\gamma(\kappa)$.

\halfrule

Next, let us revisit a construction of a metric on $M$ altered in a tubular neighbourhood around a submanifold using a curve that we call \emph{Chernysh's trick}.\footnote{While there have been similar ideas prior to the work of Chernysh (e.g.\ \cite{gajer-1987,carr-1988}), the formalism as seen in the following is essentially laid out in \cite{chernysh-2004}.}
It differs from the original construction by Gromov and Lawson in that is produces a metric on $M$ again.
\begin{optional}
  The idea is to consider curves that give rise to an embedding $M \hookrightarrow M \times \R$ instead of $M\setminus N \hookrightarrow M \times \R$.
  The pull-back metric can then be interpreted as a deformation of the metric in $M$, when previously we obtained only a complete metric on $M \setminus N$.
\end{optional}

We consider the setup as in \cref{satz:c}.
Let $g_M$ be a riemannian metric adjusted to $\phi$ on the $\overline r$-tube for an $\overline r > 0$ small enough such that $2\overline r < \operatorname{InjRad}^{\perp}_N(g_M)$, where $\operatorname{InjRad}^{\perp}_N(g_M) \coloneqq \min_{p \in N} \sup \{ r > 0 \mid \exp^{g_M, \perp}_p \text{ is injective on the } r\text{-ball}\}$ denotes the normal injectivity radius of $g_M$.
Further, let $\gamma \in \Gamma_b(\overline r)$.

Note that $g_M$ being adjusted to $\phi$ on the $\overline r$-tube implies that $\phi|_{\nu^{\leq \overline r}N}$ coincides with the normal exponential map $\exp^{\perp}$ of the metric $g_M$.

From this we can construct the embedding $\psi_{\gamma} = \psi_{\gamma}(\phi, g_M)$ given by
\begin{align*}
  \psi_{\gamma} \colon M & \to M \times \R,\\
  p & \mapsto
      \begin{cases}
        (p, 0) & \text{ if } \dop_{g_M}(p, N) \geq \overline r\\
        (\phi(r(s) \nu), t(s)) & \text{ for } p = \phi(\overline r(1 - \nicefrac{s}{b}) \nu) \text{ with } s \in (0,b), \nu \in \nu^1N\\
        (p, t(b)) & \text{ if } p \in N.
      \end{cases}
\end{align*}
whose image is denoted by
\begin{align*}
  D_{\gamma} \coloneqq \Im \psi_{\gamma} \subset M \times \R,
\end{align*}
and carries a metric $g_D$ induced from $g_M + \dop t^2$ on $M \times \R$.
Its pull-back along $\psi_{\gamma}$ to $M$ will be denoted by $g_{\gamma}$.

In fact, we obtain a continuous map
\begin{align*}
  \Gamma(\overline r) \to \Emb(M, M \times \R),
  \quad
  \gamma \mapsto \psi_{\gamma}.
\end{align*}
The embeddings in the image all coincide outside the compact set $\phi(\nu^{\leq \overline r}N)$ with the inclusion $M \hookrightarrow M \times \{0\}$.
Moreover, the embeddings $\psi_{\gamma}$ coincide for all metrics adjusted to $\phi$ on the $\overline r$-tube.

It is well-known that the pull-back of riemannian metrics along embeddings of a compact manifold is continuous with respect to the $C^{\infty}$-topology on the space of embeddings.
From this we draw the following conclusion.

Let $\{g_{\xi}\}_{\xi \in S} \subset \mathcal R(M)$ be a family of metrics on $M$, which are adjusted to $\phi$ on the $\overline r$-tube for $\overline r > 0$ such that $2\overline r < \min_{\xi \in S} \operatorname{InjRad}^{\perp}_N(g_{\xi})$.
Then we have a continuous map
\begin{align}\label{eq:deformation-by-curves}
  S \times \Gamma(\overline r) \to \mathcal R(M),
  \quad
  (\xi, \gamma) \mapsto \psi_{\gamma}^{*}(g_{\xi} + \dop t^2) =: g_{\xi, \gamma}.
\end{align}

\begin{bem}\label{bem:deforming-connection-metrics}
  Moreover, if $g_{\xi}$ is the pull-back of a connection metric on $\phi(\nu^{\leq \overline r}N)$ from $\nu N$, then $g_{\xi,\gamma}$ is a connection metric, as well.
\end{bem}

\halfrule

The curvature operator of the induced metric $g_{\gamma}$ can be connected to the curvature operator of the product metric on $M \times \R$ as demonstrated by Hoelzel.

By definition, a curve $\gamma \in \Gamma_b(\overline r)$ induces a parametrization of
$\Im \psi_{\gamma}(M \setminus N) \subset D_{\gamma}$ given by
\begin{align*}
  \gamma \colon \nu^1N \times [0,b] \to M \times \R,
  \quad
  (\nu, s) \mapsto (\phi(r(s) \nu), t(s)),
\end{align*}
at whose image points the tangent space splits as
\begin{align*}
  \Tangent_{\gamma(\nu,s)}D_{\gamma} = \Tangent_{\phi(r(s) \nu)} T(r(s)) \oplus \left<\gamma'(\nu, s)\right>,
\end{align*}
where $T(r) \coloneqq \phi(\nu^rN)$ is the distance tube around $N$.
For $(\nu_q, r) \in \nu^1N \times (0, \overline r]$, denote by $\mathcal H(\nu_q, r)$ the parallel translation of $\Tangent_qN$ into $\phi(r\nu_q)$.
Further denote by $\mathcal V(\nu_q, r)$ the orthogonal complement to $\mathcal H(\nu_q, r) \oplus \left<\partial_r\right>$, i.e.
\begin{align*}
  \mathcal H(\nu_q, r) \oplus \mathcal V(\nu_q, r) \oplus \left<\partial_r\right>
  = \Tangent_{\phi(r \nu_q)} T(r) \oplus \left<\partial_r\right>
  = \Tangent_{\phi(r \nu_q)} M.
\end{align*}
For every $(\nu, r) \in \nu^1N \times (0,\overline r]$ choose an orthonormal basis of $\mathcal H(\nu, r) \oplus \mathcal V(\nu, r)$, i.e.\ an isometry $i_{\nu, r} \colon \mathbb E^{n-1} \to \mathcal H(\nu, r) \oplus \mathcal V(\nu, r)$.\footnote{Note that in general this cannot be done a continuous way, if $\nu N$ is not assumed to be trivial.}

With respect to this choice we introduce the following notation
\begin{align}
  & \label{eq:curv-op-emb} \tilde R_D(\nu, s) \coloneqq (i_{\nu,r(s)} \oplus (-\gamma'(\nu,s)))^{*}R_D,\\
  & \nonumber \tilde R_M(\nu, r) \coloneqq (i_{\nu,r} \oplus \partial_r)^{*} R_M,\\
  & \nonumber \tilde R_T(\nu, r) \coloneqq (i_{\nu,r} \oplus \partial_t)^{*} R_T,
\end{align}
where $R_T$ is the curvature operator of $T(r) \times \R$ endowed with the product metric $g_M|_{T(r)} + \dop t^2$.

\begin{prop}[{\cite[Proposition 2.5]{hoelzel-2016}}]\label{prop:curvature-operator-on-D}
  In the situation above, for every $(\nu, s) \in \nu^1N \times I$, we have
  \begin{align*}
    \tilde R_D(\nu, s) = \cos^2 \theta(s) \tilde R_M(\nu, r(s)) + \sin^2 \theta(s) \tilde R_T(\nu, r(s)) + E(\nu, s),
  \end{align*}
  for $E(\nu, s)$ a curvature operator satisfying
  \begin{align*}
    \|E(\nu, s)\| \leq \cos \theta(s)(1 - \cos \theta(s))C_1 + \frac{\theta'(s)\sin \theta(s)}{r(s)}C_2,
  \end{align*}
  where $C_1,C_2$ are constants only depending on $D(\overline r) \coloneqq \phi(\nu^{\leq \overline r}N)$, $g_M|_{D(\overline r)}$ and $N$.
\end{prop}

If $\{g_{\xi}\}_{\xi \in S}$ is a family of metrics on $M$, we will choose isometries $i_{\nu, r}$ for every metric $g_{\xi}$.
Moreover, we denote by $\tilde R_{D,\xi}$, $\tilde R_{M,\xi}$ and $\tilde R_{T,\xi}$ the corresponding entities defined in \cref{eq:curv-op-emb}.

\subsection{Constructing the deformation map $\Pi$}

The first step will be to show that the construction of a new metric using a graph described above is actually a continuous deformation procedure within $\mathcal R_C(M)$ that can be applied to compact families of metrics simultaneously.

Throughout this section, we always consider the setup as in \cref{satz:c}.

\begin{prop}\label{prop:family-deformation}
  Let $\{g_{\xi}\}_{\xi \in S} \subset \mathcal R_C(M)$ be a family of metrics on $M$ satisfying $C$, which are adjusted to $\phi$ on the $r$-tube for some $r > 0$.
  Then there exists an $\overline r \leq r$ and a curve $\gamma \in \tilde \Gamma(\overline r)$ such that $(M,g_{\xi,\gamma})$ satisfies $C$, where $g_{\xi,\gamma}$ is obtained from $g_{\xi}$ via \cref{eq:deformation-by-curves}.

  Moreover $\gamma$ can be chosen such that according to the partition $r(s_4) = r(s_5)$ is arbitrarily small.
\end{prop}

This follows from the constructive proof of \cite[Theorem 2.1]{hoelzel-2016}, which can easily be adapted to construct a curve $\gamma$ as required for an entire family of metrics.
We only need to make sure that during the bend of $\gamma$ towards the $t$-axis $C$ remains satisfied.

\begin{lem}[{Initial bending, adapted from \cite[Lemma 2.9]{hoelzel-2016}}]\label{lemma:initial-bend}
  There exist $s_2 > 0$, $\theta_0 > 0$ and a smooth non-decreasing function $\theta \colon [0,s_2] \to [0, \theta_0]$ with $\theta'|_{[0,\varepsilon) \cup (s_2 - \varepsilon, s_2]} \equiv 0$ for all $\varepsilon > 0$ small enough such that $\tilde R_{D, \xi}(\nu, s) \in C$ for $s \in [0, s_0]$ and all $\xi \in S$.
\end{lem}

\begin{proof}
  The proof from \cite{hoelzel-2016} directly carries over to this case as all choices involved can be made in accordance with a compact family of metrics.
\end{proof}

\begin{lem}[{Second bend, adapted from \cite[Lemma 2.10]{hoelzel-2016}}]\label{lemma:second-bend}
  There exists an $r^{*} \in (0, \overline r)$ such that for every $r \in (0, r^{*})$ there is an extension of $\theta$ obtained from \cref{lemma:initial-bend} to a smooth non-decreasing function $\theta \colon [0, s_5] \to [0, \frac{\pi}{2}]$ such that $\tilde R_{D, \xi}(\nu, s) \in C$ for all $\xi \in S$, $r(s) > 0$ and $\theta|_{[s_4,s_5]} \equiv \frac{\pi}{2}$ and $r|_{[s_4,s_5]} \equiv r$ for some $s_4, s_5$ large enough.
\end{lem}

\begin{proof}
  We will only cover a part of the proof that we want to utilize later.
  Hoelzel shows in \cite[p.29f]{hoelzel-2016} that to conclude that $\tilde R_{D, \xi}(\nu, s) \in C$, it is enough to ensure
  \begin{align*}
    \theta'(s) \leq \frac{\rho}{2C_2} \frac{\sin \theta(s)}{r(s)},
  \end{align*}
  for $s \in [s_2, s_5]$, where $\rho = \rho(\tilde R_{M,\xi})$ and $C_2$ (as in \cref{prop:curvature-operator-on-D}) depend on the family of metrics.
  We let $C_3 \coloneqq \min \{ \frac{\rho}{C_2} | \, \xi \in S \}$ and conclude that while
  \begin{equation}\label{eq:second-bend}
    \theta'(s) \leq \frac{C_3}{2} \frac{\sin \theta(s)}{r(s)},
  \end{equation}
  we have $\tilde R_{D,\xi}(\theta, s) \in C$ for all $\xi \in S$.

  This can be used to inductively define the extension of $\theta$.
  Assume $\theta$ is defined on $[0, s_l]$ with $\theta(s_l) < \frac{\pi}{2}$ and define $s_{l+1} \coloneqq s_l + \frac{r(s_l)}{2}$.
  Now choose a bump function $\eta_l$ with support in $[s_l + \frac{r(s_l)}{16}, s_{l+1} - \frac{r(s_l)}{16}]$ which is constantly $\frac{C_3}{4} \frac{\sin \theta(s_l)}{r(s_l)}$ on $[s_l + \frac{r(s_l)}{8}, s_{l+1} - \frac{r(s_l)}{8}]$.
  Setting
  \begin{align*}
    \theta(s) \coloneqq \theta(s_l) + \int_{s_l}^s \eta_l(u) \dop u
  \end{align*}
  for $s \in (s_l, s_{l+1}]$ defines an extension of $\theta$ to $[0, s_{l+1}]$, which ensures that $r(s_l) \geq r(s) \geq \frac{r(s_l)}{2} > 0$ (cf. \eqref{eq:r-t-in-terms-of-theta}) and thus satisfies
  \begin{equation}\label{eq:second-bend-construction}
    \theta'(s) \leq \frac{C_3}{4} \frac{\sin \theta(s_l)}{r(s_l)}
    < \frac{C_3}{2} \frac{\sin \theta(s)}{r(s)}.
  \end{equation}
  Most importantly, $\theta$ increases at least by
  \begin{align*}
    \theta(s_{l+1}) - \theta(s_l) \geq \int_{s_l + \frac{r(s_l)}{8}}^{s_{l+1} - \frac{r(s_l)}{8}} \eta_l(u) \dop u
    \geq \frac{C_3}{16} \sin \theta_0.
  \end{align*}
  Now after finitely many steps we obtain a smooth non-increasing $\theta \colon [0,a] \to [0,\frac{\pi}{2} + \varepsilon]$, which we can adjust by a cutoff function to yield a smooth non-increasing $\theta \colon [0, a + 1] \to [0,\frac{\pi}{2}]$ with $\theta|_{[a,a+1]} \equiv 1$ keeping \eqref{eq:second-bend} satisfied.
  W.l.o.g. we can assume that $r(a) = r$, since we can let $\theta$ follow a straight line after passing $s_2$ to arbitrarily increase the $r$-coordinate.
  We let $s_4 \coloneqq a$ and $s_5 \coloneqq a + 1$.
\end{proof}

\begin{optional}
  \begin{bem*}
    Here we let $\gamma$ follow the arc of a circle such that Hoelzel Prop 2.5 has an explicit form.
  \end{bem*}
\end{optional}

\begin{proof}[Proof of \cref{prop:family-deformation}]
  By \cref{lemma:second-bend} we obtain a curve $\gamma = \gamma(\theta)$ determined by its angular function $\theta \colon [0, s_5] \to [0,\frac{\pi}{2}]$ such that $\tilde R_{D, \xi}(\nu, s) \in C$ for all $\xi \in S, \nu \in \nu^1N$ and $s \in [0,s_5]$.
  Now extend $\theta$ to $[0,s_6]$ by choosing a smooth, on $[s_5,s_6]$ non-increasing function with $\theta(s_6) = 0$ such that $\gamma(\theta)$ follows the arc of a circle (of possibly infinite radius) centered on the $t$-axis (cf. \cref{fig:bending-to-t-axis}).
  Since $\theta'|_{[s_5,s_6)} \leq 0$, \eqref{eq:second-bend}.
  Finally, we let $s_1 \coloneqq \inf\{s \geq 0 \mid \kappa(s) > 0\}$ and $s_3 \coloneqq \inf \{s \geq s_2 \mid \kappa(s) > 0 \}$ to see that $\gamma \in \tilde \Gamma(\overline r)$.
  \begin{figure}[h!]
    \centering
    \includegraphics[draft=false,width=.6\textwidth]{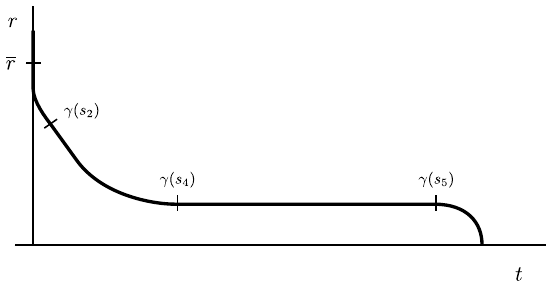}
    \caption{Bending $\gamma$ to intersect the $t$-axis following the arc of a circle.}
    \label{fig:bending-to-t-axis}
  \end{figure}
\end{proof}

\begin{prop}\label{prop:smooth-deformation}
  If $\gamma$ is obtained from \cref{prop:family-deformation}, there exists an isotopy $\alpha \colon [0,1] \to \Gamma(\overline r)$ such that $\alpha(0) = \gamma$ and $\alpha(1) = \gamma^0$, where $\gamma^0 \colon \R \to \R^2,\ s \mapsto (\overline r - s, 0)$ is the curve along the $r$-axis, such that $g_{\xi,\alpha(t)} \in \mathcal R_C(M)$ for all $\xi \in S,\ t \in [0,1]$.
\end{prop}

\begin{bem*}
  This proposition states that we can deform $\gamma$ (cf.\ \cref{fig:gl-curve}) and thereby the corresponding metrics to the originial metric keeping the curvature condition $C$ satisfied.
  Such an isotopy of curves (possibly satisfying additions assumptions) is often referred to as a \emph{Gromov-Lawson curve} (cf.\ \cite{ebert-frenck-2018,walsh-2013}).
  \begin{figure}[h!]
    \centering
    \includegraphics[draft=false,width=.6\textwidth]{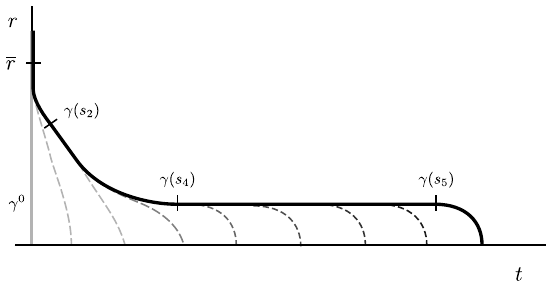}
    \caption{Isotopy between $\gamma$ and $\gamma^0$.}
    \label{fig:gl-curve}
  \end{figure}
\end{bem*}

Recall that during the so-called \emph{second bend} (within $[s_3,s_4]$), $\gamma$ satisfies \eqref{eq:second-bend} for all $s \in [s_3,s_4]$, where $C_3$ is a constant, which depends on the curvature condition $C$ and the family of metrics.
We will argue that we can modify $\gamma$ via its curvature function.
This will depend on the following Lemma adapted from \cite[Proposition 2.3]{chernysh-2004}.

\begin{lem}\label{lemma:delta-cutoff}
  Let $\gamma \in \tilde \Gamma_b(\overline r)$ and let $\kappa \colon [0,b] \to \R$ be its signed curvature function such that
  \begin{align}\label{eq:delta-cutoff-improvement}
    \theta'(s) \leq \frac{C_3}{3} \frac{\sin \theta(s)}{r(s)}
  \end{align}
  is satisfied for $s \in [s_3,s_4]$.
  Then there exists a $\delta > 0$ such that for every $s_{\bullet} \in [s_3,s_4]$ the curve $\tilde \gamma$ determined by the curvature function $\delta_{s_{\bullet}}\kappa$ satisfies \eqref{eq:second-bend}, where $\delta_{s_{\bullet}} \colon [0,b] \to [0,1]$ is a smooth $\delta$-cutoff function in the sense that
  \begin{align*}
    \delta_{s_{\bullet}}(s) =
    \begin{cases}
      1 & s \leq s_{\bullet}\\
      0 & s \geq s_{\bullet} + \delta
    \end{cases}.
  \end{align*}
\end{lem}

\begin{proof}
  The proof is entirely analogous to that of \cite[Proposition 2.3]{chernysh-2004}.
  Note that because of \cref{eq:delta-cutoff-improvement}, there exists a $\delta > 0$ such that for all $s \in [s_3, s_4]$ and $t \in [0, \delta]$ we have
  \begin{align*}
    \kappa(s + t) < \frac{C_3}{2} \frac{\sin \theta(s)}{r(s)}.
  \end{align*}
  Hence, we conclude that for all $t \in [0,\delta]$ and $s_{\bullet} \in [s_3,s_4]$
  \begin{align*}
    \tilde \kappa(s_{\bullet} + t) = \delta_{s_{\bullet}}(s_{\bullet} + t)\kappa(s_{\bullet} + t)
    \leq \kappa(s_{\bullet} + t)
    < \frac{C_3}{2} \frac{\sin \theta(s_{\bullet})}{r(s_{\bullet})}
    \leq \frac{C_3}{2} \frac{\sin \tilde \theta(s_{\bullet} + t)}{\tilde r(s_{\bullet} + t)}.
  \end{align*}
  For $s \in [s_3, s_{\bullet}]$ both $\kappa$ and $\tilde \kappa$ coincide, while for $t > \delta$, we have $\tilde \kappa(s_{\bullet} + t) = 0$.
  Thus, \eqref{eq:second-bend} is satified.
\end{proof}

\begin{proof}[Proof of \cref{prop:smooth-deformation}]
  We will argue in two steps.
  First we will show that we can deform a curve $\gamma$ that bends up to a straight line of small angle $\theta_0$ within $[s_2,s_3]$ and bends down to meet the $t$-axis in a right angle to $\gamma^0$ maintaining $C$.

  Using \cref{prop:curvature-operator-on-D} Hoelzel concludes that for all $\nu \in \nu^1N$ and $s \in [0,b)$
  \begin{align*}
    \|\tilde R_D(\nu, s) - \tilde R_M(\nu, r(s))\|
    & \leq \sin^2 \theta(s)(\sup \|R_M\| + \sup \|R_T\|)|_{\phi(r(s)\nu)}\\
    & \quad + \cos \theta(s)(1-\cos\theta(s))C_1 + \frac{\theta'(s) \sin\theta(s)}{r(s)}C_2,
  \end{align*}
  where the suprema are taken over the points in the tubular neighbourhood $D(\overline r)$.
  As this is true for all $s \in [0,s_2]$ (where $\theta(s) \leq \theta_0$), it is easy to see that it remains satisfied, if we linearly decrease $\theta$ to $0$.

  Let $p \colon [0,1] \to I,\ t \mapsto s_4 - t(s_4 - s_3)$ be the linear path from $s_4$ to $s_3$.
  Choose $\delta > 0$ from \cref{lemma:delta-cutoff} (whose assumptions are satisfied, which can be seen from \eqref{eq:second-bend-construction}) and a $\delta$-cutoff function $\delta_{s_4} \colon I \to [0,1]$.
  Note that we obtain a continuous family of $\delta$-cutoff functions $\tilde \delta_{s_{\bullet}}(s) \coloneqq \delta_{s_4}(s + (s_4 - s_{\bullet}))$ depending on $s_{\bullet} \in [s_3,s_4]$.
  Now define for $t \in [0,1]$
  \begin{align*}
    \kappa_t(s) \coloneqq
    \begin{cases}
      \delta_{p(t)}(s) \kappa(s) & \text{ if } 0 \leq s \leq p(t) + \delta\\
      0 & \text{ if } p(t) + \delta \leq s \leq s_t\\
      \delta_{s_t}(s_t + \delta - s)\varepsilon_t \kappa(s + (s_5 - s_t)) & \text{ if } 0 \leq s_t \leq L_t
    \end{cases}
  \end{align*}
  where $s_t$ is the intersection of $\gamma(\tilde \delta_{p(t)}\kappa(s))$ with $\{r_4\} \times \R$ and $\varepsilon_t, L_t$ are uniquely determined such that
  \begin{align*}
    \int_{[s_t, L_t]} \kappa_t(s) \dop s = - \int_{[0,s_t]} \kappa_t(s) \dop s
    \quad \text{ and } \quad
    (\gamma(\kappa_t))(L_t) \in \{0\} \times \R.
  \end{align*}
  The resulting isotopy $t \mapsto \gamma(\kappa_t)$ deforms $\gamma = \gamma(\kappa_0)$ into a curve of the form discussed in the beginnning.
\end{proof}

\begin{kor}
  There exists a curve  $\gamma \in \tilde \Gamma(\overline r)$ and a continuous map
  \begin{align*}
    \mathfrak A \colon S \times [0,1] \to \mathcal R_C(M),
  \end{align*}
  with $\mathfrak A(\xi, 0) = g_{\xi}$ and $\mathfrak A(\xi,1) = g_{\xi,\gamma}$.
\end{kor}

\begin{proof}
  By \cref{prop:smooth-deformation}, we obtain an isotopy of curves and thus an isotopy of embeddings $M \hookrightarrow M \times \R$ whose corresponding metrics obtained via \cref{eq:deformation-by-curves} satisfy $C$.
\end{proof}

In the next step, we will show that we can deform metrics to become rotationally symmetric in a small normal tube.

\begin{prop}\label{prop:middle-stage-deformation}
  There exists an $r_{*} > 0$  and a continuous map
  \begin{align*}
    \mathfrak B \colon S \times [0,1] \to \mathcal R_C(M)
  \end{align*}
  such that $\mathfrak B(\xi, 0) = g_{\xi}$ and $(\phi|_{\nu^{\leq r_{*}}N})^{*}\mathfrak B(\xi,1)$ is the restriction of a connection metric on $\nu N$.

  Moreover, if $g_{\xi}$ is contained in $\mathcal R_C^{\rot}(M)$, then $t \mapsto \mathfrak B(\xi, t)$ is a path within $\mathcal R_C^{\rot}(M)$.
\end{prop}

During the proof we will utilize the following technical result by Hoelzel.

\begin{satz}[{\cite[Proposition 2.7]{hoelzel-2016}}]\label{satz:small-tube}
  Let $C$ be be a codimension $c$ surgery stable curvature condition and let $N^k \subset M^n$ be a compact submanifold of the riemannian manifold $(M,g_M)$ with codimension $n - k \geq c$.
  There exists an $r_{*} > 0$ such that for all $r \in (0,r_{*})$ the riemannian manifold $(T(r) \times \R, g_{T(r)} + g_{\R})$ satisfies $C$.
  Moreover, there exists an $L > 0$ such that
  \begin{align*}
    \tilde R_T(\nu, r) = B_{\nicefrac{L}{r}}(\nicefrac{1}{r^2}R_{\E^{k+1} \times \mathbb S^{n-k-1}}) \subset C.
  \end{align*}
\end{satz}

We note that, since all the metrics $g_{\xi}$ are adjusted to $\phi$ on the $\overline r$-tube, we have $T(r) = \phi(\nu^rN)$ for every $r \leq \overline r$ when we apply this theorem to $g_{\xi}$.

\begin{lem}\label{lemma:gamma-circle}
  Let $\gamma = (\delta \cos s, \delta \sin s)$ for a $s \in [0,\nicefrac{\pi}{2})$ be the curve in the $(r,t)$-plane following the arc of a circle.
  Let $g_M$ be a riemannian metric on $M$ and let $D_{\gamma} \subset \phi(\nu^{\leq \delta}N) \times \R \subset M \times \R$ be obtained from $\gamma$.
  Then there exists a $\delta > 0$, such that $\tilde R_D(\nu, s)$ satisfies $C$ for all $\nu \in \nu^1N$ and $s \in [0,\nicefrac{\pi}{2})$.
\end{lem}

\begin{proof}
  The angular function $\theta$ is given by $\cos \theta(s) = \delta \sin(s)$.
  By \cref{prop:curvature-operator-on-D}, we have that
  \begin{align*}
    \tilde R_D(\nu, s) & = \cos^2\theta \tilde R_M(\nu, r(s)) + \sin^2 \theta \tilde R_T(\nu, r(s)) + E(\nu, s)\\
                       & = \delta^2\sin^2s \tilde R_M(\nu, r(s)) + (1 - \delta^2\sin^2s) \tilde R_T(\nu, r(s)) + E(\nu, s),
  \end{align*}
  where $\|E(\nu, s)\| \leq \delta \sin s(1 - \delta \sin s)C_1 - \frac{\delta \cos s}{\sqrt{1-\delta^2\sin^2s}} C_2$.
  Moreover, we conclude from \cref{satz:small-tube} that for $\delta$ small enough there exists an $L > 0$ such that
  \begin{align*}
    \tilde R_T(\nu, r(s))
    & = \frac{1}{r(s)^2} R_{\E^{k+1} \times \mathbb S^{n-k-1}} + \tilde E(\nu, s)\\
    & = \frac{1}{\delta^2\cos^2s} R_{\E^{k+1} \times \mathbb S^{n-k-1}} + \tilde E(\nu, s)
  \end{align*}
  where $\|\tilde E(\nu,s)\| \leq \frac{L}{r(s)} = \frac{L}{\delta \cos s}$.
  Combined, we have
  \begin{align*}
    \tilde R_D(\nu, s) & = \delta^2\sin^2s \tilde R_M(\nu, r(s)) + \frac{1 - \delta^2\sin^2s}{\delta^2\cos^2s} R_{\E^{k+1} \times \mathbb S^{n-k-1}}\\
                       & \qquad + (1-\delta^2\sin^2s) \tilde E(\nu, s) + E(\nu, s).
  \end{align*}
  If we choose $\delta$ small enough we ensure that $\lambda(\delta) \coloneqq \frac{1 - \delta^2\sin^2s}{\delta^2\cos^2s}$ is large enough such that $\lambda(\delta) R_{\E^{k+1} \times \mathbb S^{n-k-1}}$ is contained in $C$.
  Because $\|(1-\delta^2\sin^2s) \tilde E(\nu, s)\|$ grows slower than the cone opens, for a $\delta$ small enough $\tilde R_D(\nu, s)$ is contained in $C$.
\end{proof}

\begin{proof}[Proof of \cref{prop:middle-stage-deformation}]
  Note that $\gamma \in \tilde \Gamma(\overline r)$ as obtained by \cref{prop:family-deformation} comes with a partition $0 = s_0 \leq s_1 \leq \cdots \leq s_6 = b$.
  We denote $t_i = t(s_i)$ and define $I_1 \coloneqq (-\infty, t_4]$, $I_2 = [t_4, t_5]$ and $I_3 = [t_5, \infty)$.
  As before, let $D(\overline r) \coloneqq \phi(\nu^{\leq \overline r}N)$ and define $S_j \coloneqq D(\overline r) \times I_J$.
  This partitions $D_{\gamma}$ into $D_j \coloneqq D_{\gamma} \cap S_j$ for $j = 1,2,3$ with the following descriptions (cf.\ \cref{fig:partition})
  \begin{align*}
    D_1 \coloneqq & \{ (\phi(r(s)\nu), t(s)) \mid 0 \leq s \leq s_4, \nu \in \nu^1 N \}, \\
    D_2 \coloneqq & \{ (\phi(r(s)\nu), t(s)) \mid s_4 \leq s \leq s_5, \nu \in \nu^1 N \},\\
    D_3 \coloneqq & \{ (\phi(r(s)\nu), t(s)) \mid s_5 \leq s < b, \nu \in \nu^1 N \} \cup \{ (p, t(b)) \mid p \in N \}.
  \end{align*}

  Let $h^{\torp}$ denote a connection metric on $\nu N$ obtained from a torpedo metric (as constructed in \cref{sec:connection-metrics}).

  \begin{figure}[h!]
    \centering
    \includegraphics[draft=false,width=.8\textwidth]{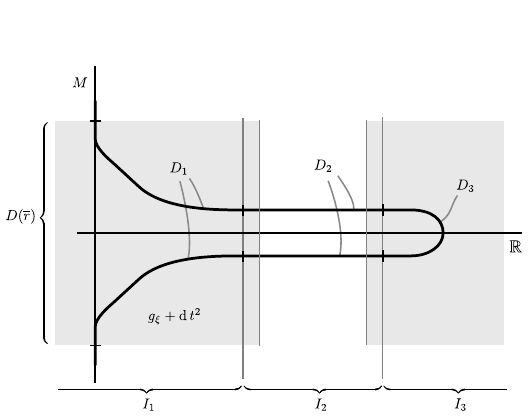}
    \caption{Partition of $D_{\gamma}$ into $D_1, D_2$ and $D_3$.}
    \label{fig:partition}
  \end{figure}

  Let $\xi \in S$ and let $h$ be the pull-back of a connection metric on $\nu N$, which is adjusted to $\phi$ on the $r$-tube, along $(\phi|_{\nu^{\leq \overline r}N})^{-1}$ to $\phi(\nu^{\leq \overline r}N)$.

  Fix a small $\varepsilon > 0$ and define for $l \in [0,1]$ a riemannian metric
  \begin{align*}
    G_{\xi, l} \coloneqq
    \begin{cases}
      g_{\xi} + \dop t^2 & \text{ on } D(\overline r) \times (I_1 \cup [t_4, t_4 + \varepsilon))\\
      ((1-l)g_{\xi} + lh) + \dop t^2 & \text{ on } D(\overline r) \times ((t_5 - \varepsilon, t_5] \cup I_3).
    \end{cases}
  \end{align*}
  For every $l \in [0,1]$ continue this to a smooth metric on $D(\overline r) \times \R$ by choosing a family of smooth paths of metrics $P_l \colon [t_4+\varepsilon, t_5-\varepsilon] \to \mathcal R(M)$, which are adjusted to $\phi$ on some tube.
  By an application of \cref{prop:adjusting-to-phi}, we can assume that the paths are through metrics adjusted to $\phi$ on the $r_0$-tube for some $r_0 < \overline r$.
  Moreover, we choose a path through connection metrics, in case $g_{\xi}$ is a connection metric to begin with.

  In total, we obtain a family of metrics $\{G_{\xi,l}\}_{(\xi, l) \in S \times [0,1]}$.

  Since the metrics $G_{\xi,l}$ restricted to $M \times \{t\}$ are adjusted to $\phi$ on the $r_0$-tube, $D_2$ is the distance tube around $S_2 \cap (N \times \{0\})$ for every $(\xi,l) \in S \times [0,1]$.
  Thus, by \cref{satz:small-tube} the metric restricted to $D_2$ satisfies $C$, if $r(s_4) = r(s_5) =: \delta$ is small enough.
  This can be accomplished according to the construction of $\gamma$ (cf. \cref{prop:family-deformation}).

  In $S_3$, the submanifold $D_3$ is determined by $\gamma$ following the arc of a circle.
  Hence, we can apply \cref{lemma:gamma-circle} to see that for $\delta$ small enough $C$ is satified for all $(\xi, l) \in S \times [0,1]$.
  Because $h$ is a connection metric and modified by $\gamma$ following an arc only radially, the metric induced on $D_3$ by $G_{\xi, 1}$ is a connection metric, as well.

  Therefore, the metric induced on $D_{\gamma}$ from $g_{\xi, l}$ satisfies $C$ and if we let $r_{*} \coloneqq \overline r(1 - \nicefrac{t_4}{b})$, then
  \begin{align*}
    \mathfrak B \colon S \times [0,1] \to \mathcal R_C(M),
    \quad
    (\xi, l) \mapsto
    \begin{cases}
      \mathfrak A(\xi, 2l) & \text{ for } l \in [0,\nicefrac{1}{2}], \\
      \phi_{\gamma}^{*}(G_{\xi, 2l-1}|_{D_{\gamma}}) & \text{ for } l \in [\nicefrac{1}{2},1]
    \end{cases}
  \end{align*}
  has the properties claimed.

  In particular, if $g_{\xi} \in \mathcal R_C^{\rot}(M)$, then by \cref{bem:deforming-connection-metrics} the metrics $\mathfrak A(\xi, \blank)$ are contained in $\mathcal R_C^{\rot}(M)$.
  During the second part of the deformation the restriction of $G_{\xi,l}$ to $M \times \{t\}$ is a connection metric and thus $\mathfrak B(\xi, \blank)$ is contained in $\mathcal R_C^{\rot}(M)$.
\end{proof}

The final step amounts to an adjustment of the metric we produced to the tubular map $\phi$ using a suitable radial diffeomorphism.

\begin{proof}[Proof of \cref{prop:deformation-of-families}]
  Consider the deformation map $\mathfrak B$ as constructed in \cref{prop:middle-stage-deformation}.
  We can choose a family of radial diffeomorphisms $\Phi \colon [0,1] \to \Diff(M),\ t \mapsto \Phi_t$, such that
  \begin{enumerate}
  \item $\Phi_0 \equiv \id_M$,
  \item $\Phi_t$ is the identity outside of $\phi(\nu^{\leq \nicefrac{3}{2}\overline r}N)$,
  \item $\Phi_1(\phi(\nu^{\leq \overline r}N)) \subset \phi( \nu^{\leq r_{*}}N)$.
  \end{enumerate}
  Now define the final deformation map
  \begin{align*}
    \Pi \colon [0,1] \times S \to \mathcal R_C(M),
    \quad
    (l, \xi) \mapsto
    \begin{cases}
      \mathfrak B(2l, \xi) & \text{ for } l \in [0,\nicefrac{1}{2}] \\
      \Phi_{2l-1}^{*}\mathfrak B(1, \xi) & \text{ for } l \in [\nicefrac{1}{2},1]
    \end{cases}.
  \end{align*}
  This map has the desired properties and in particular $\Pi(1, \xi) = \Phi_1^{*}\mathfrak B(1, \xi)$ is contained in $\mathcal R_C^{\rot}(M)$.
\end{proof}

\section{Rotationally symmetric metrics}\label{section:rotationally-sym-metrics}

The main goal of this section is the proof of \cref{prop:deformation-warped-to-torpedo}, which will complete our proof of the main theorem.
It will follow from an analysis of rotationally symmetric metrics on the disc.
As before, we are generalizing the method laid out in Chernysh's \cite{chernysh-2004} to the case of certain curvature conditions.

\subsection{Metrics on the disc}

In the following fix a deformable, codimension $c$ surgery stable curvature condition $C \subset \mathcal C_{\mathrm B}(\E^n)$, $\delta > 0$ and $q \geq c$.

\begin{dfn}\label{dfn:rc-rot}
  Define
  \begin{align*}
    \mathcal R^{\rot} \coloneqq \{ g \in \mathcal R(\R^q) \mid g = \alpha^2(t)\dop t^2 + \beta^2(t) g_{\mathbb S^{q-1}} \text{ and } \alpha(t) \neq 0 \ \forall t, \beta \geq 0 \}
  \end{align*}
  in spherical coordinates $(0, \infty] \times S^{q-1} \cong S^{q}\setminus\{0\}$, where $\alpha,\beta \colon (0,\delta] \to \R$.
  Further, we denote by
  \begin{align*}
    \mathcal R^{\rot}_C \coloneqq \{ g \in \mathcal R^{\rot} \mid g_{\mathbb E^{n-q}} + g \in \mathcal R_C(\R^{n-q} \times \R^q) \}
  \end{align*}
  the space of rotationally symmetric metrics satisfying $C$.
\end{dfn}

\begin{dfn}
  On $\mathcal R^{\rot}$ the \emph{radius} defines a continuous map
  \begin{align*}
    \operatorname{rad} \colon \mathcal R^{\rot} \to \R_{> 0},
    \quad
    g \mapsto \int_0^{\delta} g(\gamma', \gamma')^{\frac{1}{2}}\dop t,
  \end{align*}
  where $\gamma$ is a radial path from the centre to a boundary point of $D^q(\delta)$.
  In the spherical coordinates $(0, \operatorname{rad}(g)] \times S^{n-1} \cong D^n(\delta)\setminus\{0\}$ the metric $g$ is of the form $\dop t^2 + \beta g_{\mathbb S^{n-1}}$ for $\beta \colon (0, \operatorname{rad} g] \to \R$, which is the restriction of a smooth odd function $\beta \colon \R \to \R$ with $\beta'(0) = 1$ and $\beta^{(\text{even})}(0) = 0$.
  The group of diffeomorphisms $\Diff_0(\R) = \{ \psi \in \Diff(\R) \mid \psi^{(k)}(0) = 0 \ \forall k \geq 0, \psi'(0) = 1\}$ acts continuously on rotationally symmetric metrics as follows:
  \begin{align*}
    \Diff_0(\R) \times \mathcal R^{\rot} & \to \mathcal R^{\rot},\\
    (\psi, g = \dop t^2 + \beta^2g_{\mathbb S^{q-1}}) & \mapsto \psi \star g \coloneqq \dop t^2 + (\beta \circ \psi)^2 g_{\mathbb S^{q-1}},
  \end{align*}
  where the expression for $\psi \star g$ is understood in spherical co-ordinates $(0, \psi^{-1}(\operatorname{rad}(g))] \times S^{q-1} \cong D^n(\delta)\setminus\{0\}$.
\end{dfn}

\begin{bem}
  \begin{enumerate}
  \item For the radius we have $\operatorname{rad}(\psi \star g) = \psi^{-1}(\operatorname{rad} g)$.
  \item If $R^{\beta}$ is the curvature operator corresponding to $g = \dop t^2 + \beta^2g_{\mathbb S^{q-q}}$ (cf.\ \eqref{eq:curvature-operator-beta}), then $R^{\beta \circ \psi}$ is the curvature operator corresponding to $\psi \star g$.
  \end{enumerate}
\end{bem}

\begin{prop}\label{prop:psi-1}
  In the situation of the above definition there exists a continuous function $\sigma \colon \mathcal R^{\rot}_C \to (0, \frac{\delta}{2}]$ and a continuous deformation of the identity $\Psi_1 \colon \mathcal R^{\rot}_C \times [0,1] \to \mathcal R^{\rot}_C$ with
  \begin{enumerate}
  \item $\Psi_1(\blank, 0) \equiv \id$,
  \item $\Psi_1(\blank, s) \equiv \id$ near $\partial D^q(\delta)$ for all $s \in [0,1]$,
  \item $\Psi_1(g, 1) = \dop t^2 + \beta(t) g_{\mathbb S^{q-1}}$ with
    \begin{enumerate}[label=(\roman*)]
    \item $\beta^{(l)}|_{\sigma(g)} = 0$ for all $l \geq 1$,
    \item $0 \leq \beta' \leq 1$ and $\beta'' \leq 0$ on $[0, \sigma(g)]$.
    \end{enumerate}
  \end{enumerate}
\end{prop}

The proof will be based on the observation that the warping function of a metric $g \in \mathcal R^{\rot}_C$ behaves in a certain way near $0$ combined with explicit deformations of the metric corresponding to this region.
We will separate these steps in the following lemmata.

\begin{lem}\label{lemma:tip}
  For every $g = \dop t^2 + \beta^2g_{\mathbb S^{q-1}} \in \mathcal R^{\rot}_C$ there exists a $t^{*}$ such that $0 < \beta' < 1$ and $\beta'' < 0$ on $(0, t^{*}]$.
\end{lem}

\begin{proof}
  Assume this would not be the case, i.e.\ there exists a point $t_0$ with $\beta' \geq 1$ or $\beta'' \geq 0$ on $[0,t_0]$.
  Since we have $\beta'(0) = 1$ and $\beta''(0) = 0$, both inequalities are actually equivalent to convexity of $\beta$ on the interval $[0,t_0]$.
  The curvature operator corresponding to $\beta$ is
  \begin{equation}\label{eq:warped-operator}
    R^{\beta} = \underbrace{\frac{1-\beta'^2}{\beta^2}}_{=:\lambda(\beta)} R_{\mathbb E \times \mathbb S^{q-1}} \underbrace{- \frac{\beta''}{\beta}}_{=: \mu(\beta)} L_q
  \end{equation}
  and for $t \in [0, t_0]$ we thus have $\lambda(\beta)(t) \leq 0$ and $\mu(\beta)(t) \leq 0$.

  Since we know that $C$ satisfies an inner ray condition with respect to $s(R_{\mathbb E \times \mathbb S^{q-1}} + \mu L)$ for $\mu \geq 0$ and $s > 0$, we conclude that with $\lambda \leq 0,\ \mu \leq 0$ the curvature condition $C$ would contain the flat curvature operator $0$, which is a contradiction to our assumption that $C$ is deformable.
\end{proof}

\begin{lem}[{\cite[Lemma 3.5]{chernysh-2004}}]\label{lemma:first-deform}
  For $0 < C_1 \leq 1$, $0 < t^{*}$, $0 < t_{*} < \nicefrac{1}{2}\ t^{*}$, $s \in [0,1]$ let $0 < a = a(t_{*}, t^{*}) < b = b(t_{*}, t^{*}) < t_{*}$ be continuous functions.
  There exist continuous functions $0 < C_2 = C_2(C_1, t^{*}) \leq 1$, $e = e(s)$ and an isotopy through diffeomorphisms $\phi_s \colon \R \to \R$ such that for $b < c = \phi_s^{-1}(\nicefrac{8}{10}\ t^{*}) < d = \phi_s^{-1}(\nicefrac{9}{10}\ t^{*}) < e$
  \begin{enumerate}[label=(\roman*)]
  \item $\phi_0 \equiv \id$,
  \item $\phi_s(0) = 0$ and $\phi_s(e(s)) = t^{*}$,
  \item $\phi_s'|_{[0,a] \cup [d,\infty)} \equiv 1$ and $\phi_s'|_{[b,c]} \equiv 1 - sC_2$, $\phi_s^{(n)}|_{[b,c]} \equiv 0$ for all $n \geq 2$,
  \item $0 \leq 1 - C_2 \leq \phi_s' \leq 1$,
  \item $\phi_s'' \leq C_1$,
  \item $\phi_s''|_{[a,b]} \leq 0$ and $\phi_s''|_{[c,d]} \geq 0$.
  \end{enumerate}
\end{lem}

\begin{figure}[h!]
  \centering
  \includegraphics[draft=false]{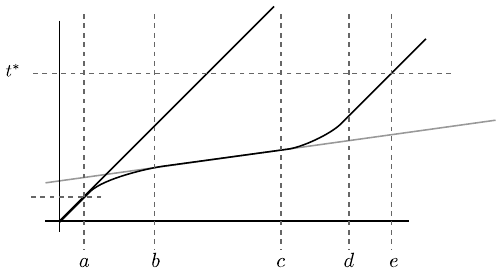}
  \caption{First deformation family as described in \cref{lemma:first-deform}}
\end{figure}

\begin{lem}[{\cite[Lemma 3.7]{chernysh-2004}}]\label{lemma:second-deform}
  Let $0 < D_1 \leq 1$, $0 < t^{**}$, $s \in [0,1]$.
  There exist continuous functions $0 < \overline b = \overline b(D_1, t^{**}) < \frac{t^{**}}{2}$, $\overline e = \overline e(s)$ and an isotopy through diffeomorphisms $\psi_s \colon \R \to \R$ such that for $0 < \overline a = \nicefrac{1}{10}\ \overline b < \overline c = \psi_s^{-1}(\nicefrac{9}{10}\ t^{**}) < \overline e$
  \begin{enumerate}[label=(\roman*)]
  \item $\psi_0 \equiv \id$,
  \item $\psi_s(0) = 0$, $\psi_s(\overline d(s)) = t^{**}$,
  \item $0 \leq \psi_s' \leq 1$,
  \item $\psi_s'|_{_{[0,\overline a]\cup[\overline c,\infty)}} \equiv 1$,
  \item $\psi_s''|_{(0, \infty)}(t) \leq \nicefrac{D_1}{t}$,
  \item $\psi_1^{(n)}(\overline b(1)) = 0$ for all $n \geq 1$.
  \end{enumerate}
\end{lem}

\begin{figure}[h!]
  \centering
  \includegraphics[draft=false]{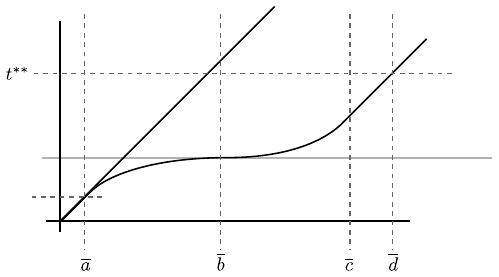}
  \caption{Second deformation family as described in \cref{lemma:second-deform}}
\end{figure}

The proof of these two lemmata is a straight forward explicit construction by means of ordinary differential equations.

\begin{proof}[Proof of \cref{prop:psi-1}]
  We consider $g = \dop t^2 + \beta^2g_{\mathbb S^{q-1}} \in \mathcal R^{\rot}_C$.
  By \cref{lemma:tip}, there exists a $t^{*}$ such that $0 < \beta' < 1$ and $\beta'' < 0$ on $(0, t^{*}]$.
  Without loss of generality, we can assume that $t^{*} \leq \nicefrac{\delta}{2}$.

  Now we invoke \cref{lemma:first-deform} with the parameters $C_1$, $t^{*}$ and the \cref{lemma:second-deform} with the parameters $D_1, t^{**}$ to define
  \begin{align*}
    \Psi_1 \colon (g, s) \mapsto
    \begin{cases}
      \phi_{2s} \star \beta & \text{ for } s \in [0, \nicefrac{1}{2}]\\
      (\phi_1 \circ \psi_{2s-1}) \star \beta & \text{ for } s \in [\nicefrac{1}{2},1]\\
    \end{cases}.
  \end{align*}

  We need to make sure that during both parts of the deformation the corresponding metrics stay within the space $\mathcal R_C^{\rot}$.
  \begin{enumerate}
  \item We need to make sure that the rotationally symmetric metric given by the warping function $\beta \circ \phi_{2s}$, which has the curvature operator
    \begin{equation}\label{eq:deformed-curvature-operator}
      R^{\beta \circ \phi_{2s}} = \frac{1 - \beta'(\phi_{2s})^2 \phi_{2s}'^2 }{\beta(\phi_{2s})^2} R_{\mathbb E \times \mathbb S^{q-1}} - \frac{\beta''(\phi_{2s}) \phi_{2s}'^2 + \beta'(\phi_{2s}) \phi_{2s}''}{\beta(\phi_{2s})} L_q,
    \end{equation}
    satisfies $C$, as well.
    On $[0,a]$, we have $\phi_{2s} \equiv \id$ and there is nothing to check.
    On $[d,\phi_{2s}^{-1}(\operatorname{rad} g)]$ we have $\phi_{2s}' \equiv 1$ and $\phi_{2s}'' \equiv 0$, thus, $R^{\beta \circ \phi_{2s}}$ at $t$ equals $R^{\beta}$ at $\phi_{2s}^{-1}(t)$, which of course satisfies $C$.

    \noindent
    On $[b,c]$ we have
    \begin{align*}
      R^{\beta \circ \phi_{2s}} & = \frac{1-\beta'(\phi_{2s})^2(1-sC_2)^2}{\beta(\phi_{2s})^2} R_{\mathbb E \times \mathbb S^{q-1}} - \frac{\beta''(\phi_{2s})(1 - sC_2)^2}{\beta(\phi_{2s})}L_q\\
                                & = (1 - sC_2)^2 R^{\beta}|_{\psi_{2s}} \in C.
    \end{align*}
    In the remaining part $[a,b] \cup [c,d]$ we will use the deformability (cf. \cref{dfn:deformability}) of the curvature condition $C$.
    We can conclude that the first summand in \eqref{eq:deformed-curvature-operator} is contained in $C$ because $\frac{1 - \beta'^2(\phi_{2s}) \phi_{2s}'^2 }{\beta^2(\phi_{2s})} > 0$.
    The whole sum is clearly contained in $C$, if $ - \frac{\beta''(\phi_{2s}) \phi_{2s}'^2}{\beta(\phi_{2s})} - \frac{\beta'(\phi_{2s}) \phi_{2s}''}{\beta(\phi_{2s})} \geq 0$, i.e.\ if $-\frac{\beta''(\phi_{2s})\phi_{2s}'^2}{\beta'(\phi_{2s})} \geq \phi_{2s}''$.
    Since $\beta'' < 0$, $1 - C_2(C_1, t^{*}) \leq \phi_{2s}' \leq 1$ and $0 < \beta' < 1$, we can choose $C_1$ in \cref{lemma:first-deform} such that this is satisfied.
  \item Now we need to conclude the same for $\beta \circ \phi_1 \circ \psi_{2s-1}$, where this time we denote $\tilde \beta \coloneqq \beta \circ \phi_1$.
    We let $b < t^{**} < c$ from \cref{lemma:first-deform}.
    The curvature operator here is
    \begin{align}
      \nonumber R^{\beta \circ \phi_1 \circ \psi_{2s-1}} &  = \frac{1 - (\tilde \beta')^2(\psi_{2s-1}) (\psi_{2s-1}')^2}{\tilde \beta^2(\psi_{2s-1})} R_{\mathbb E \times \mathbb S^{q-1}}\\
      \label{eq:deformed-curvature-operator-2}  & \qquad - \frac{\tilde\beta''(\psi_{2s-1}) \psi_{2s-1}'^2 + \tilde \beta'(\psi_{2s-1}) \psi_{2s-1}''}{\tilde \beta(\psi_{2s-1})} L_q
    \end{align}
    again, by the same reasoning, the first summand is contained in $C$.
    And again on $[0,\overline a] \cup [\overline c,(\phi_1 \circ \psi_{2s-1})^{-1}(\operatorname{rad} g)]$ there is nothing to check.
    Within $[\overline a, \overline c]$ by deformability, the sum in \eqref{eq:deformed-curvature-operator-2} would be contained in $C$, if $-\frac{\tilde \beta'' \psi_{2s-1}'^2}{\tilde \beta} \geq \psi_{2s-1}''$.
    Unfortunately, the left-hand side of this inequality is not necessarily strictly positive and hence can't be used to saturate $D_1$ in \cref{lemma:second-deform}.
    Using the fact that $C$ satisfies an inner cone condition with respect to $R_{\mathbb E \times \mathbb S^{q-1}}$, we conclude that actually the factor in front of $L_q$ has to be larger than a negative constant depending on the cone opening.

    \begin{figure}[h!]
      \centering
      \includegraphics[draft=false]{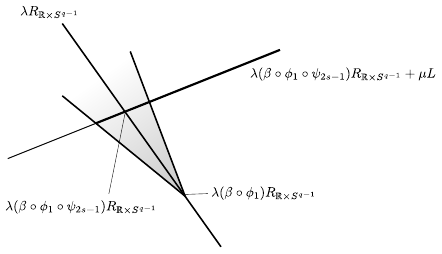}
      \caption{Intersection of the cone condition with the $\lambda R_{\mathbb E \times \mathbb S^{q-1}} + \mu L_q$ plane. Here $\lambda(\beta) \coloneqq \frac{1-\beta'^2}{\beta^2}$.}
      \label{fig:cone-opening}
    \end{figure}

    On $[\overline a,\overline c]$ we have $ \tilde \beta'' = \beta''(\phi_1) \phi_1'^2 + \beta'(\phi_1) \phi_1'' < 0$, because $\phi_1'' \leq 0$, $\beta'' < 0$, $0 < \phi_1' \leq 1$ and $0 < \beta' < 1$, as well as $0 < \tilde \beta' = \beta'(\phi_1) \phi_1' < 1$.
    Now \cref{eq:deformed-curvature-operator-2} is contained in $C$, if
    \begin{align*}
      -\frac{\tilde \beta''(\psi_{2s-1})\psi_{2s-1}'^2 + \tilde \beta'(\psi_{2s-1})\psi_{2s-1}''}{\tilde \beta} \geq -C^{**},
    \end{align*}
    where $C^{**} > 0$ depends on the cone opening $\rho(\lambda(\beta \circ \phi_1) R_{\mathbb E \times \mathbb S^{q-1}}) > 0$ (cf. \cref{dfn:inner-cone-condition} and \cref{fig:cone-opening}) and the difference $\lambda(\beta \circ \phi_1 \circ \psi_{2s-1}) - \lambda(\beta \circ \phi_1) > 0$.
    This is equivalent to
    \begin{align*}
      \frac{\tilde \beta(\psi_{2s-1})C^{**} - \tilde \beta''(\psi_{2s-1}) \psi_{2s-1}'^2}{\tilde \beta'} \geq \psi_{2s-1}'',
    \end{align*}
    where the left-hand side is strictly positive, which enables us to choose $D_1$ such that $R^{\beta \circ \phi_1 \circ \psi_{2s-1}}$ is contained in $C$.
  \end{enumerate}
  Finally, let $\sigma(g) \coloneqq b$ from \cref{lemma:second-deform}.
\end{proof}

\begin{dfn}\label{dfn:rc-loc}
  Let $g_{\torp} = \dop t^2 + \beta_{\torp}(t)^2g_{\mathbb S^{n-k-1}}$ be the torpedo metric on $D^{n-k}(\delta)$.
  We define
  \begin{align*}
    \textstyle \mathcal R_C^{\operatorname{loc}} \coloneqq \{ g \in \mathcal R_C^{\rot} \mid g = \big(\frac{\delta}{\delta^{*}}\big)^2 \dop t^2 + \beta_{\torp}\big(\frac{\delta}{\delta^{*}}t\big)^2 g_{\mathbb S^{q-1}} \text{ with } \delta^{*} \in (0,\delta] \}
  \end{align*}
  and refer to it as the \emph{space of locally torpedo metrics}.
  The radius $\delta^{*}$ on which a metric $g \in \mathcal R_C^{\operatorname{loc}}$ is torpedo yields a continuous map $\delta^{*} \colon \mathcal R_C^{\operatorname{loc}} \to (0,\delta]$.
\end{dfn}

\begin{prop}\label{prop:rot-to-loc}
  In the situation of \cref{prop:psi-1} there exists a continuous map
  \begin{align*}
    \Psi_2 \colon \mathcal R^{\rot}_C \times [0,1] \to \mathcal R^{\rot}_C
  \end{align*}
  such that
  \begin{enumerate}
  \item $\Psi_2(\blank, 0) \equiv \id$,
  \item $\Psi_2(\blank, s) \equiv \id$ near $\partial D^q(\delta)$ for all $s \in [0,1]$,
  \item $\Psi_2(g, 1) = \left(\frac{\delta}{\sigma}\right)^2 \dop t^2 + \beta_{\torp}\left(\frac{\delta}{\sigma}t\right)^2 g_{\mathbb S^{q-1}}$ on $D^q(\sigma)$,
  \end{enumerate}
\end{prop}

\begin{proof}
  For the first part of the deformation we use $\Psi_1$, i.e.
  \begin{align*}
    \Psi_2 \colon (g, s) \mapsto \Psi_1(g, 2s) & \text{ for } s \in [0,\nicefrac{1}{2}].
  \end{align*}
  Now let $\Psi_1(g, 1) = \dop t^2 + \beta^2g_{\mathbb S^{q-1}}$.
  We will describe a deformation on $D^q(\sigma(g))$.
  For $u \in [0,1]$ let $\beta_u \coloneqq (1-u) \beta + u \beta_{\torp}$ we observe the following.
  If $\beta(t) \leq \beta_u(t) \leq \beta_{\torp}(t)$, we have
  \begin{align*}
    \frac{1 - \beta_{\torp}'^2(t)}{\beta_{\torp}(t)}
    \leq \frac{1 - \beta_u'^2(t)}{\beta_u(t)}
    \leq \frac{1 - \beta'^2(t)}{\beta(t)},
  \end{align*}
  since $0 \leq \beta' \leq 1$ and $0 \leq \beta_{\torp}' \leq 1$ (we obtain the opposite inequalities if $\beta(t) \geq \beta_u(t) \geq \beta_{\torp}(t)$).
  Moreover
  \begin{align*}
    \frac{\beta_u''}{\beta_u} = \frac{(1-u) \beta + u \beta_{\torp}}{(1-u) \beta + u \beta_{\torp}} \leq 0,
  \end{align*}
  because $\beta'' \leq 0$ and $\beta_{\torp} \leq 0$.
  Hence, by deformability $\beta_u$ is contained in $C$ and since we have created a collar at $\sigma(g)$ (i.e.\ $\beta^{(l)}(\sigma(g)) = 0$ for all $l \geq 1$), we can deform $\beta$ on a small interval to connect $\beta_u(\sigma(g))$ with $\beta(\sigma(g))$.
  Denote the result by $\tilde \beta_u$ and define $\Psi_2 \colon (g,s) \mapsto \dop t^2 + \tilde \beta_{2s-1}^2 g_{\mathbb S^{q-1}}$ for $s \in [\nicefrac{1}{2},1]$.
\end{proof}

\subsection{Rotationally symmetric metrics around a submanifold}

Finally, we return to the case of an embedded submanifold with the assumptions as listed in \cref{satz:c}.

Recall that we can use $g_N, h^{\nu N}, \omega$ and a rotationally symmetric metric $g_{\rot}$ on $\R^{n-k}$ to define a connection metric $h^{\nabla}(g_{\rot}) := g_N \oplus_\omega g_{\rot}$ on $\nu N$.

We note that every rotationally symmetric metric around $N$ in the sense of \cref{dfn:metrics-in-standard-form} is given by a rotationally symmetric on the disc, that is contained in $\mathcal R_C^{\rot}$, in particular.
More precisely, for $\delta = R$ in \cref{dfn:rc-rot} we have for $g \in R_C^{\rot}(M)$ that $\phi|_{\nu^{\leq R}N}^{*}g = h^{\nabla}(g_{\rot})|_{\nu^{\leq R}N}$, where $g_{\rot}|_{D^{n-k}(R)} \in \mathcal R_C^{\rot}$.
As $g$ is adjusted to $\phi$ on the $R$-tube, we have that $g_{\rot} = \dop r^2 + \beta(r)^2g_{\mathbb S^{n-k-1}}$.

We note that a continuous alteration of $g_{\rot}$, which does not alter the metric near $S^{n-k}(R) \subset \mathcal R^{n-k}$ results in a continuous alteration of the metric $g$ near $N$.

\begin{dfn}
  Analogous to \cref{dfn:rc-loc} we define
  \begin{align*}
    \mathcal R_C^{\operatorname{loc}}(M) \coloneqq \{ g \in \mathcal R_C(M) \mid \phi|_{\nu^{\leq R}N}^{*}g = h^{\nabla}(g_{\rot})|_{\nu^{\leq R}N} \text{ with } g_{\rot}|_{D^{n-k}(R)} \in \mathcal R_C^{\operatorname{loc}} \}.
  \end{align*}
\end{dfn}

\begin{lem}\label{lemma:loc-to-torp}
  There exists a continuous deformation map $\Phi \colon \mathcal R_C^{\operatorname{loc}}(M) \times [0,1] \to \mathcal R_C^{\operatorname{loc}}(M)$ such that
  \begin{enumerate}[label=(\roman*)]
  \item $\Phi(\blank, 0) \equiv \id$,
  \item $\Phi(\blank, 1) \in \mathcal R_C^{\torp}(M)$.
  \end{enumerate}
\end{lem}

\begin{proof}
  For a small $\varepsilon > 0$ choose a family of radial diffeomorphisms on the disc $\phi \colon D^{n-k}(R + \varepsilon) \times (0, R] \times [0,1]  \to D^{n-k}(R + \varepsilon)$ such that
  \begin{enumerate}[label=(\roman*)]
  \item $\phi(\blank, r, 0) \equiv \id$,
  \item $\phi(x, r, 1) = \frac{r}{R}x$ on $x \in D^{n-k}(R)$.
  \item Denote by $\tilde \delta = \tilde \delta(r)$ the value where $\phi(\blank, r, 1)^{-1}(S^{n-k}(R)) = S^{n-k}(\tilde \delta)$.
    If $\tilde \delta \in [\delta, \delta + \varepsilon]$, then $\phi(\blank, r, 1)$ is a radial euclidean isometry in a tubular neighbourhood of $S^{n-k}(\tilde \delta)$.
  \end{enumerate}
  These give rise to a family of diffeomorphisms on the normal bundle $\phi \colon \nu^{\leq R + \varepsilon}N \times (0,R] \times [0,1] \to \nu^{\leq R + \varepsilon}N$.
  Further, by \cref{dfn:rc-loc} we obtain a continuous function $\delta^{*} \colon \mathcal R_C^{\operatorname{loc}}(M) \to (0, \delta]$ such that we can define
  \begin{align*}
    \Phi(g, s) \coloneqq \phi(\blank, \delta^{*}(g), s)^{*}g,
  \end{align*}
  which is the deformation map we were seeking to construct.
\end{proof}

Now we are set up to give a proof of \cref{prop:deformation-warped-to-torpedo}, i.e.\ to show that $\mathcal R_C^{\torp}(M)$ is a weak deformation retract of $\mathcal R_C^{\rot}(M)$.

\begin{proof}[Proof of \cref{prop:deformation-warped-to-torpedo}]
  The idea is to construct a retraction map
  \begin{align*}
    r \colon \mathcal R_C^{\rot}(M) \xrightarrow{\Psi} \mathcal R_C^{\operatorname{loc}}(M) \xrightarrow{\Phi} \mathcal R_C^{\torp}(M),
  \end{align*}
  where $\Psi$ is a replacement of the rotationally symmetric metric by local torpedo metric, while $\Phi$ is induced by a family of radial diffeomorphisms on the disc.

  The map $\Psi_2$ constructed in \cref{prop:rot-to-loc} gives rise to a continuous deformation map $\Psi \colon \mathcal R_C^{\rot}(M) \times [0,1] \to \mathcal R_C^{\rot}(M)$ with $\Psi(\blank,1) \subset \mathcal R_C^{\operatorname{loc}}(M)$ by replacing the rotationally symmetric metric in the connection metric with the metrics obtained during the deformation $\Psi_2$.
  Now use this $\Psi$ and the map $\Phi$ constructed in \cref{lemma:loc-to-torp} to define $r \coloneqq \Phi(\blank, 1) \circ \Psi_2(\blank, 1)$.
  The homotopy $i \circ r \simeq \id_{\mathcal R_C^{\rot}(M)}$, where $i$ is the inclusion, is given by
  \begin{align*}
    (g, s)
    \mapsto
    \begin{cases}
      \Psi(g, 2s) & \text{ for } s \in [0, \nicefrac{1}{2}]\\
      \Phi(\Psi(g,1), 2s-1) & \text{ for } s \in [\nicefrac{1}{2}, 1]
    \end{cases}.
  \end{align*}
  For the other homotopy $r \circ i \simeq \id_{\mathcal R_C^{\torp}(M)}$, we note that $\Psi_2(\blank, 1)$ does not alter the torpedo metric within the annulus $D^{n-k}(\nicefrac{R}{2}, R)$, since $t^{*} \leq \nicefrac{R}{2}$ in the proof of \cref{prop:psi-1}.
  Let $g \in \mathcal R_C^{\torp}(M)$.
  On $D^{n-k}(R)$ the metric $r(g)$ is torpedo, i.e.\ agrees with $\tau^{*}g$ , while on the annulus $D^{n-k}(R, R + \varepsilon)$ the metric might not agree with the original one.
  This can be resolved by straightening this to a collar region and shrinking it down afterwards, which is possible because $\tilde \delta \in [R, R + \varepsilon]$ in the proof of \cref{lemma:loc-to-torp}.
\end{proof}

\begin{acknowledgements}
  I would like to express my gratitude to Prof.\ W.\ Tuschmann for his constant support during the preparation of this work, which is based on the results of my Ph.D.\ thesis.
  Moreover, I would like to thank Prof.\ J.\ Ebert,  Dr.\ M.\ Wiemeler and Dr.\ G.\ Frenck for helpful discussions on the subject.
  In particular, I would like to thank Dr.\ M.\ Wiemeler for his suggestion to consider Chernysh's construction with non-trivial normal bundles.
\end{acknowledgements}

%%% Local Variables:
%%% mode: latex
%%% TeX-master: "paper-homotopy-invariance-of-R_C-arxiv"
%%% End:

\bibliographystyle{amsalpha}
\bibliography{literature}

\end{document}